\documentclass[final,leqno,onefignum,onetabnum]{siamltex1213}

\title{SISO Output Affine Feedback Transformation \\ Group and Its Fa\`{a} di Bruno Hopf Algebra}

\author{W.~Steven Gray\thanks{Old Dominion University, Norfolk, Virginia 23529, USA ({\tt sgray@odu.edu}).}
\and Kurusch Ebrahimi-Fard\thanks{Norwegian University of Science and Technology -- NTNU, 7491 Trondheim, Norway
         ({\tt kurusch.ebrahimi-fard@math.ntnu.no}). On leave from UHA, Mulhouse, France.}}

\usepackage{cite} 
\usepackage{color} 
\usepackage{graphicx} 
\usepackage{latexsym,amssymb,amsmath} 
\usepackage{mathrsfs} 
\usepackage{dsfont} 
\usepackage{cancel} 

\allowdisplaybreaks[4]

\def\abs#1{\left\vert #1 \right\vert}
\def\allpoly{\mbox{$\re\langle X \rangle$}}

\def\allpolyx0degn{\mbox{$P_n$}}

\def\allseries{\mbox{$\re\langle\langle X \rangle\rangle$}}
\def\allseriesdelta{\mbox{$\re\langle\langle X_\delta \rangle\rangle$}}
\def\allseriesdeltaLC{\mbox{$\re_{LC}\langle\langle X_\delta \rangle\rangle$}}

\def\allseriesdeltanp{\mbox{$\re_{\rm np}\langle\langle X_\delta \rangle\rangle$}}
\def\allseriesell{\mbox{$\re^{\ell} \langle\langle X \rangle\rangle$}}
\def\allserieselldelta{\mbox{$\re^\ell\langle\langle X_\delta \rangle\rangle$}}

\def\allseriesLC{\mbox{$\re_{LC}\langle\langle X \rangle\rangle$}}

\def\allseriesmLC{\mbox{$\re^{m}_{LC}\langle\langle X \rangle\rangle$}}

\def\allseriesnLC{\mbox{$\re^{n}_{LC}\langle\langle X \rangle\rangle$}}

\def\allseriesnp{\mbox{$\re_{\rm np}\langle\langle X \rangle\rangle$}}

\def\allseriesellLC{\mbox{$\re^{\ell}_{LC}\langle\langle X \rangle\rangle$}}

\def\allseriesXOm{\mbox{$\re^m [[ X_0 ]]$}} 
\def\allseriesXOmLC{\mbox{$\re^m_{LC} [[ X_0 ]]$}}

\def\bfem#1{{\bf \em #1}} 

\newcommand{\comment}[1]{} 

\def\dim{{\rm dim}}
\def\dist{{\rm dist}}

\def\emptyword{\emptyset}

\def\Endallseries{{\rm End}(\allseries)}

\def\eqref#1{(\ref{#1})} 

\def\Fliessdelta{\mathscr{F}_{\delta}}

\def\id{{\rm id}}

\def\liepoly{{\cal L}(X)}

\def\mbf#1{\hbox{\mathversion{bold}$#1$}} 
\def\modcomp{\:\tilde{\circ}\,} 

\def\norm#1{\left\Vert#1\right\Vert}
\def\notin{{\not\in}}

\def\ord{{\rm ord}}

\def\re{{\mathbb R}} 

\def\sameau{\rule[0.017in]{0.2in}{0.012in}}

\def\shuffle{{\scriptscriptstyle \;\sqcup \hspace*{-0.05cm}\sqcup\;}}

\def\supp{{\rm supp}}

\newtheorem{example}[theorem]{Example}

\def\begce{\begin{center}}
\def\endce{\end{center}}
\def\begar{\begin{array}}
\def\endar{\end{array}}
\def\begeq{\begin{equation}}
\def\endeq{\end{equation}}
\def\begdi{\begin{displaymath}}
\def\enddi{\end{displaymath}}
\def\begdis{\begin{eqnarray*}}
\def\enddis{\end{eqnarray*}}
\def\begeqa{\begin{eqnarray}}
\def\endeqa{\end{eqnarray}}
\def\begdes{\begin{description}}
\def\enddes{\end{description}}
\def\begit{\begin{itemize}}
\def\endit{\end{itemize}}
\def\begen{\begin{enumerate}}
\def\enden{\end{enumerate}}
\def\beglar{\left[\begin{array}}
\def\endrar{\end{array}\right]}
\def\begle{\begin{lemma}}
\def\endle{\end{lemma}}
\def\begde{\begin{definition}}
\def\endde{\end{definition}}
\def\begth{\begin{theorem}}
\def\endth{\end{theorem}}
\def\begco{\begin{corollary}}
\def\endco{\end{corollary}}
\def\begprop{\begin{proposition}}
\def\endprop{\end{proposition}}
\def\begex{\begin{example}}
\def\endex{\end{example}}
\def\begexer{\begin{exercise}}
\def\endexer{\end{exercise}}

\def\begres{\noindent{\bf Remarks}:\begin{enumerate}}
\def\endres{\end{enumerate} \par}
\def\begpr{\begin{proof}}
\def\endpr{\end{proof}}

\def\begprx#1{\par{\it Proof #1}. \ignorespaces} 
\def\endprx{\vbox{\hrule height0.6pt\hbox{%
   \vrule height1.3ex width0.6pt\hskip0.8ex
   \vrule width0.6pt}\hrule height0.6pt
  }}
\def\begtab{\begin{tabular}}
\def\endtab{\end{tabular}}
\def\rref#1{(\ref{#1})}

\begin{document}

\maketitle
\slugger{sicon}{xxxx}{xx}{x}{x--x} 

\begin{abstract}
The general goal of this paper is to identify a transformation group that can be used
to describe a class of feedback interconnections involving subsystems which are modeled solely in terms of
{\em Chen--Fliess functional expansions} or {\em Fliess operators} and are independent of the existence of any state space
models.
This interconnection, called an {\em output affine feedback connection}, is distinguished from conventional output feedback
by the presence of a multiplier in an outer loop.
Once this transformation group is established,
three basic questions are addressed. How can this transformation group be used to provide an explicit Fliess operator representation
of such a closed-loop system? Is it possible to use this feedback scheme to do system inversion purely in an input-output setting?
In particular, can feedback input-output linearization be posed and solved entirely in this framework, i.e., without
the need for any state space realization? Lastly, what can be
said about feedback invariants under this transformation group?
A final objective of the paper is to describe the Lie algebra of infinitesimal characters associated with
the group in terms of a pre-Lie product.
\end{abstract}

\begin{keywords}
nonlinear control systems, Chen--Fliess series, Hopf algebras, transformation groups
\end{keywords}

\begin{AMS}
93C10, 93B18, 16T30
\end{AMS}

\pagestyle{myheadings}
\thispagestyle{plain}
\markboth{W. STEVEN GRAY AND KURUSCH EBRAHIMI-FARD}{SISO OUTPUT AFFINE FEEDBACK TRANSFORMATION GROUP}

\section{Introduction}

Feedback transformation groups have been used extensively in control theory
since its inception to explain the way that feedback can alter the nature of a system.
The early work of Brockett and Krishnaprasad in
the case of linear systems \cite{Brockett_83,Brockett-Krishnaprasad_80}
and Brockett, Jakubczyk, Respondek, and many others in the context of nonlinear state space
systems \cite{Brockett_78,Jakubczyk_90,Jakubczyk-Respondek_80,Tall-Respondek_03} has been
important, for example, in identifying feedback invariants.
The general goal of this paper is to identify a transformation group that can be used
to describe a class of feedback interconnections involving subsystems which are modeled solely in terms of
{\em Chen--Fliess functional expansions} or {\em Fliess operators} and are independent of the existence of any state space
models \cite{Fliess_81,Fliess_83}.
\begin{figure}[t]
\begin{center}
\includegraphics[scale=0.26]{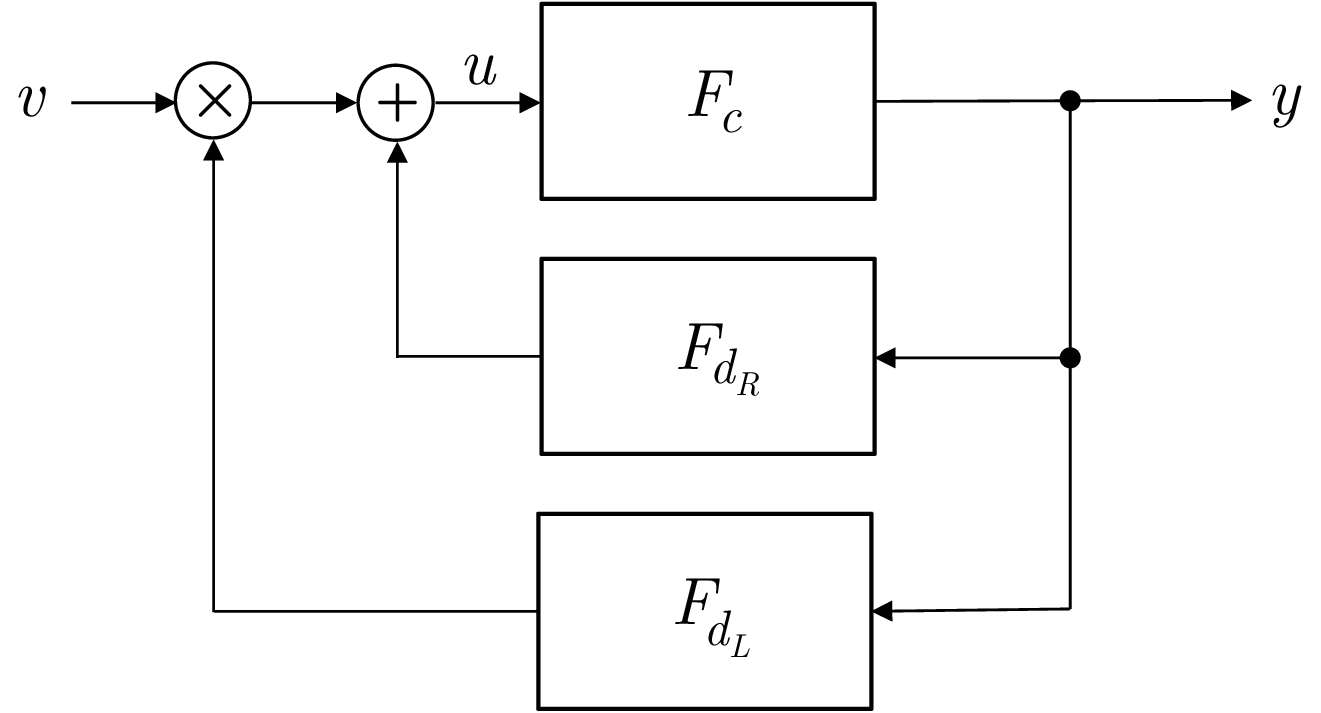}
\caption{Output affine feedback connection of Fliess operators}
\label{fig:affine-output-feedback-connection}
\end{center}
\end{figure}
This interconnection, called an {\em output affine feedback connection}, is shown in Figure~\ref{fig:affine-output-feedback-connection}.
The presence of a multiplier in an outer loop distinguishes this structure from conventional output feedback.
Once this transformation group is established,
three basic questions are addressed. How can this transformation group be used to provide an explicit Fliess operator representation
of such a closed-loop system? Is it possible to use this feedback scheme to do system inversion purely in an input-output setting?
In particular, can feedback input-output linearization be posed and solved entirely in this framework, i.e., without
the need for any state space realization? Lastly, what can be
said about feedback invariants under this transformation group?
A final objective of the paper is to describe the Lie algebra of infinitesimal characters associated with
this new transformation group and show that it is induced by a pre-Lie
product analogous to what was found for the output feedback transformation group in \cite{Foissy_CiA15,Foissy_EJM15}.
This fact may give deeper insight into the combinatorial aspects of the theory and group invariants (see \cite{Gray-etal_CDC14}),
but this aspect of the problem will be deferred to a later publication.

The basic approach to solving these problems is to use an operator algebra encoded by various products of
formal power series.
This is somewhat analogous to the familiar way in which interconnections of
linear time-invariant systems are analyzed using only transfer functions and the
algebra of rational functions. The feedback connection in general requires the use of a certain inverse
operator that is naturally described in terms of a group.
To do any calculations explicitly and
efficiently here, it is necessary to present an associated Hopf algebra induced by the group product and its inverse.
It has been known for a long time in some physics communities that Hopf algebras provide
a natural computational framework for working with groups \cite{Cartier_07,Figueroa-Gracia-Bondia_05,Frabetti-Manchon_15}. For example, the group inverse
can be computed explicitly in terms of  the {\em antipode}, a particular antihomomorphism of the Hopf algebra. If the latter is both connected and graded, then this antipode has a recursive description which is ideal for software implementation.
As will be seen, the Hopf algebra of interest here is directly related to a composition of Fliess operators.
The inversion formula that it provides is a kind of functional analogue of Lagrange's formula for the inversion of smooth
functions in terms of their Taylor series coefficients. This classical formula has a combinatorial interpretation in terms of a connected graded Hopf algebra \cite{Joni-Rota_79}. In particular, the coproduct of this Hopf algebra is defined in terms of
the classical Fa\`{a} di Bruno formula which gives the Taylor series coefficients for the composition of two sufficiently well behaved functions.
Therefore, this Hopf algebra is often referred to as a {\em Fa\`{a} di Bruno type} Hopf algebra \cite{Figueroa-Gracia-Bondia_05,Frabetti-Manchon_15}. The Hopf algebra presented in this paper is also of Fa\`{a} di Bruno type and
is best viewed as a generalization of existing work on output feedback (corresponding to setting $F_{d_L}[y]=1$ for all
$y$ in Figure~\ref{fig:affine-output-feedback-connection}) by the authors and others appearing in \cite{Duffaut-etal_JAC16,Foissy_CiA15,Foissy_EJM15,Gray-Duffaut_Espinosa_CDC13,Gray-etal_Auto14,Gray-etal_CDC14,Gray-etal_SCL14,Gray-Wang_08,Thitsa-Gray_12}.
The focus here will be on the single-input single-output (SISO) case as the multivariable case introduces additional
technical issues that are nontrivial and require more advanced concepts from Hopf algebra theory.
Finally, it should be noted that a portion of this paper appeared in preliminary form in \cite{Gray_MTNS14}.

The paper is organized as follows. In the next section, some preliminary concepts concerning Fliess operators, ultrametric spaces,
and Hopf algebras are briefly outlined.
In Section~\ref{sec:transformation-group-allseriesdelta},
the output affine feedback transformation group is described. The problem of determining a Fliess operator representation of the
output affine closed-loop system is then solved.
The Hopf algebra is next developed in Section~\ref{sec:hopf-algebra}
and subsequently applied in the following section
to address the remaining two problems concerning system inversion/feedback linearization and feedback invariants.
In Section~\ref{sec:FPS-Lie-group}, the Lie algebra of infinitesimal characters is described.

\section{Preliminaries}
\label{sec:preliminaries}

A finite nonempty set of noncommuting symbols $X=\{ x_0,x_1,$ $\ldots,x_m\}$ is
called an {\em alphabet}. Each element of $X$ is called a {\em
letter}, and any finite sequence of letters from $X$,
$\eta=x_{i_1}\cdots x_{i_k}$, is called a {\em word} over $X$. The
{\em length} of the word $\eta$, $\abs{\eta}$, is given by the number of letters it
contains.
The set of all words with length $k$ is denoted by
$X^k$. The set of all words including the empty word, $\emptyset$,
is designated by $X^\ast$, while $X^+:=X^\ast-\{\emptyset\}$. The set $X^\ast$
forms a monoid under catenation.
The set $\eta X^\ast$ is comprised of all words with the prefix $\eta\in X^\ast$.
For any fixed integer $\ell\geq 1$, a
mapping $c:X^\ast\rightarrow
\re^\ell$ is called a {\em formal power series}. The value of $c$ at
$\eta\in X^\ast$ is written as $(c,\eta)\in\re^{\ell}$ and called the {\em coefficient} of
$\eta$ in $c$.
Typically, $c$ is
represented as the formal sum $c=\sum_{\eta\in X^\ast}(c,\eta)\eta.$
If the {\em constant term} $(c,\emptyset)=0$ then $c$ is said to be {\em proper}.
The {\em support} of $c$, $\supp(c)$, is the set of all words having nonzero coefficients.
The {\em order} of $c$, $\ord(c)$, is the length of the shortest word in its support ($\ord(0):=\infty$).\footnote{For notational convenience,
$p=(p,\emptyset)\emptyset\in\allpoly$ is often abbreviated as $p=(p,\emptyset)$.}
The collection of all formal power series over $X^*$ is denoted by
$\allseriesell$.
It forms an associative $\re$-algebra under
the catenation product and a commutative and associative $\re$-algebra under the shuffle
product, denoted here by the shuffle symbol $\shuffle$.
The latter is the
$\re$-bilinear extension of the shuffle product of two words,
which is defined inductively by
\begin{equation}
\label{shuffleproduct}
(x_i\eta)\shuffle
(x_j\xi)=x_i(\eta\shuffle(x_j\xi))+x_j((x_i\eta)\shuffle \xi)
\end{equation}
with
$\eta\shuffle\emptyset=\emptyset\shuffle\eta=\eta$ for all $\eta,\xi\in X^\ast$ and $x_i,x_j\in X$
\cite{Fliess_81}.

\subsection{Fliess Operators and Their Interconnections}

One can formally associate with any series $c\in\allseriesell$ a causal
$m$-input, $\ell$-output operator, $F_c$, in the following manner.
Let $\mathfrak{p}\ge 1$ and $t_0 < t_1$ be given. For a Lebesgue measurable
function $u: [t_0,t_1] \rightarrow\re^m$, define
$\norm{u}_{\mathfrak{p}}=\max\{\norm{u_i}_{\mathfrak{p}}: \ 1\le
i\le m\}$, where $\norm{u_i}_{\mathfrak{p}}$ is the usual
$L_{\mathfrak{p}}$-norm for a measurable real-valued function,
$u_i$, defined on $[t_0,t_1]$.  Let $L^m_{\mathfrak{p}}[t_0,t_1]$
denote the set of all measurable functions defined on $[t_0,t_1]$
having a finite $\norm{\cdot}_{\mathfrak{p}}$ norm and
$B_{\mathfrak{p}}^m(R)[t_0,t_1]:=\{u\in
L_{\mathfrak{p}}^m[t_0,t_1]:\norm{u}_{\mathfrak{p}}\leq R\}$. Assume $C[t_0,t_1]$
is the subset of continuous functions in $L_{1}^m[t_0,t_1]$. Define
inductively for each word $\eta\in X^{\ast}$ the map $E_\eta:
L_1^m[t_0, t_1]\rightarrow C[t_0, t_1]$ by setting
$E_\emptyset[u]=1$ and letting
\[E_{x_i\bar{\eta}}[u](t,t_0) =
\int_{t_0}^tu_{i}(\tau)E_{\bar{\eta}}[u](\tau,t_0)\,d\tau, \] where
$x_i\in X$, $\bar{\eta}\in X^{\ast}$, and $u_0=1$. The
input-output operator corresponding to $c=\sum_{\eta\in X^\ast}(c,\eta)\eta\in\allseriesell$ is the {\em Fliess operator}
\begdi
F_c[u](t) =
\sum_{\eta\in X^{\ast}} (c,\eta)\,E_\eta[u](t,t_0) \label{eq:Fliess-operator-defined}
\enddi
\hspace*{-0.05in}\cite{Fliess_81,Fliess_83}.
If there exist real numbers $K_c,M_c>0$ such that
\begeq
\abs{(c,\eta)}\le K_c M_c^{|\eta|}|\eta|!,\;\; \forall\eta\in X^{\ast},
\label{eq:local-convergence-growth-bound}
\endeq
then $F_c$ constitutes a well defined mapping from
$B_{\mathfrak p}^m(R)[t_0,$ $t_0+T]$ into $B_{\mathfrak
q}^{\ell}(S)[t_0, \, t_0+T]$ for sufficiently small $R,T >0$ and some $S>0$,
where the numbers $\mathfrak{p},\mathfrak{q}\in[1,\infty]$ are
conjugate exponents, i.e., $1/\mathfrak{p}+1/\mathfrak{q}=1$  \cite{Gray-Wang_02}.
(Here, $\abs{z}:=\max_i \abs{z_i}$ when $z\in\re^\ell$.) The set of all such
{\em locally convergent} series is denoted by $\allseriesellLC$.
A Fliess operator $F_c$ defined on $B_{\mathfrak p}^m(R)[t_0,t_0+T]$
is said to be {\em realizable} when there exists
a state space model consisting of $n$ differential
equations and $\ell$ output functions
\begin{subequations}
\begin{align}
\dot{z}&= g_0(z)+\sum_{i=1}^m g_i(z)\,u_i,\;\;z(t_0)=z_0  \label{eq:state}\\
y_j&=h_j(z),\;\;j=1,2,\ldots,\ell, \label{eq:output}
\end{align}
\end{subequations}
where each $g_i$ is an analytic vector field expressed in local
coordinates on some neighborhood ${\cal W}$ of $z_0$,
and each output function $h_j$ is an analytic function on ${\cal W}$ such that
(\ref{eq:state}) has a well defined solution $z(t)$, $t\in[t_0,t_0+T]$ on ${\cal W}$
for any given input $u\in B_{\mathfrak
p}^m(R)[t_0,t_0+T]$, and $y_j(t)=F_{c_j}[u](t)=h_j(z(t))$, $t\in[t_0,t_0+T]$, $j=1,2,\ldots,\ell$
\cite{Fliess_81,Fliess_83,Gray-Wang_02,Isidori_95,Nijmeijer-vanderSchaft_90}.
It can be shown that for any word $\eta=x_{i_k}\cdots x_{i_1}\in X^\ast$
\begeq
(c_j,\eta)=L_{g_{\eta}}h_j(z_0):=L_{g_{i_1}}\cdots L_{g_{i_k}}h_j(z_0), \label{eq:c-equals-Lgh}
\endeq
where $L_{g_i}h_j$ is the {\em Lie derivative} of $h_j$ with respect to $g_i$.
For any $c\in\allseriesell$, the $\re$-linear mapping
${\cal H}_c:\allpoly\rightarrow \allseriesell$
uniquely specified by
$
({\cal H}_c(\eta),\xi)=(c,\xi\eta),\;\; \xi,\eta \in X^{\ast}
$
is called the {\em Hankel mapping} of $c$. The series $c$ is said to
have {\em Lie rank} $n$ when the range space of ${\cal H}_c$
restricted to the $\re$-vector space of Lie polynomials over $X$,
i.e., the free Lie algebra $\liepoly$, has dimension $n$.
It is well known that $F_c$ is realizable if and
only if $c\in\allseriesellLC$ has finite Lie rank
\cite{Fliess_83,Isidori_95,Nijmeijer-vanderSchaft_90}.

When convergent Fliess operators are connected in a parallel-sum configuration, it is elementary
to show that the composite system satisfies $F_c+F_d=F_{c+d}$.
It was
shown in \cite{Fliess_81} for the parallel-product connection that $F_cF_d=F_{c\shuffle d}$.
If $F_c$ and $F_d$ with $c\in\allseriesellLC$ and
$d\in\allseriesmLC$ are interconnected in a cascade manner, the composite
system $F_c\circ F_d$ has the
Fliess operator representation $F_{c\circ d}$, where
$c\circ d$ denotes
the {\em composition product} of $c$ and $d$ as described in
\cite{Ferfera_80}.
This product is associative and $\re$-linear in its left argument $c$.
In the event that two convergent Fliess operators are interconnected
to form a feedback system,
the closed-loop system has a convergent Fliess operator representation whose generating
series is the {\em output feedback product} of $c$ and $d$, denoted by $c@d$ \cite{Gray-etal_SCL14,Thitsa-Gray_12}.
Consider, for example, the SISO case where the alphabet $X=\{x_0,x_1\}$ and $\ell=1$. Define the set of {\it unital Fliess operators}
$
	\Fliessdelta=\{I+F_c:c\in\allseriesLC\},
$
where $I$ denotes the identity map, i.e., $I[u]=u$.
It is convenient to introduce the symbol
$\delta$ as the (fictitious) generating series for $I$. That is,
$F_\delta:=I$ such that $I+F_c:=F_{\delta+c}=F_{c_\delta}$ with
$c_\delta:=\delta+c$.
The set of all such generating series for
$\Fliessdelta$ is denoted by $\allseriesdeltaLC$ with the augmented alphabet $X_{\delta}:=\{\delta\} \cup X$. It has been proved that
$\Fliessdelta$ forms a group
under the composition $F_{c_\delta}\circ F_{d_\delta}=(I+F_c)\circ(I+F_d)= F_{c_\delta\circ d_\delta}$,
where $c_\delta\circ d_\delta:=\delta+d+c\modcomp d$,
and $\modcomp$ denotes the {\em modified composition product} \cite{Gray-Li_05}.\footnote{The same symbol will be
used for composition on $\allseries$ and $\allseriesdelta$.
As elements in these two sets have a distinct
notation, i.e., $c$ versus $c_\delta$, respectively, it will always be clear which product is at play.}
The unit element of this group is the identity map $I$.
Of central importance is the so-called {\it{coordinate algebra}} of functions, $H$, on the corresponding
group $(\allseriesdeltaLC,\circ,\delta)$. In this setting the group multiplication is dualized to a comultiplication on $H$, which is captured through the notion of a Hopf algebra. In fact, the compositional nature of the product on $\Fliessdelta$ turns  $H$ into a Fa\`{a} di Bruno type Hopf algebra. In which case, the group (composition) inverse in $(\allseriesdeltaLC,\circ,\delta)$, $c_\delta^{-1}=:\delta+c^{-1}$, can be computed efficiently by a recursive algorithm \cite{Duffaut-etal_JAC16}. This inverse is also key in describing the output feedback product as shown in the following theorem.

\begth {\rm{\cite{Gray-etal_SCL14,Thitsa-Gray_12}}}
\label{th:feedback-product-formula}
For any $c,d\in\allseriesdeltaLC$ it follows that
$
c@d=c\modcomp(-d\circ c)^{-1} \label{eq:catd-formula}
$
with $c@d\in\allseriesdeltaLC$. In addition, $c@d=c\modcomp(d\circ(c@d))$.
\endth

It is known that the set of locally convergent generating series is closed under addition, the shuffle product \cite{Wang_90},
the composition product \cite{Gray-Li_05}, the modified composition product \cite{Li_04}, and the feedback product \cite{Gray-etal_SCL14,Thitsa-Gray_12}.
Therefore, all the corresponding elementary interconnections of Fliess
operators preserve local convergence.

Finally, the algebraic theory described above remains largely intact for formal power series
that do not necessarily satisfy convergence condition \rref{eq:local-convergence-growth-bound}. In place of a convergent Fliess operator
one can introduce a {\em formal} Fliess operator which maps formal functions to formal functions \cite{Gray-Wang_08}.
The motivation for this is the observation that if $c\in\allseriesellLC$ and $u$ is analytic with a convergent Taylor series
representation $u=\sum_{n\geq 0} (c_u,x_0^n)(t-t_0)^n/n!=\sum_{n\geq 0}(c_u,x_0^n)E_{x_0^n}[u](t,t_0)$, then $y=F_c[u]$ is also analytic with
Taylor series coefficients given by $c_y=c\circ c_u$ so that $c_u\in\allseriesXOmLC$ and $c_y \in \mathbb{R}_{LC}^\ell[[X_0]]$ with $X_0:=\{x_0\}.$
But when $c$ is not locally convergent,
the mapping $\allseriesXOm\rightarrow\mathbb{R}_{LC}^\ell[[X_0]]:c_u\mapsto c_y=c\circ c_u$
is still well defined. It can be viewed as the formal version of $F_c$.  The set $\allserieselldelta$, which also forms a group
under the product $c_\delta\circ d_\delta$, is the corresponding set of generating series for these formal operators.

\subsection{Shuffle Product Operations on Ultrametric Spaces}
Given a set $S$, a function ${\mathscr U}:S\times S\rightarrow \re$
is called an {\em ultrametric} if it satisfies the following
properties for all $s,s^{\prime},s^{\prime\prime}\in S$:
\begin{description}
\item[i.] ${\mathscr U}(s,s^{\prime})\geq 0$
\item[ii.] ${\mathscr U}(s,s^{\prime})=0$ if and only if $s=s^{\prime}$
\item[iii.] ${\mathscr U}(s,s^{\prime})={\cal {\mathscr U}}(s^{\prime},s)$
\item[iv.] ${\mathscr U}(s,s^{\prime})\leq \max\{{\mathscr U}(s,s^{\prime\prime}),{\mathscr U}(s^{\prime\prime},s^{\prime})\}$.
\end{description}
The pair $(S,{\mathscr U})$ is referred to as an {\em ultrametric space}.
In the event that property {\rm iv} is replaced with the triangle
inequality,
$(S,{\mathscr U})$ reduces to the usual definition of a metric space.
It is not difficult to show for any fixed $0<\sigma<1$ that
the $\re$-vector space $\allseries$ with the distance between two series defined as
$\dist(c,d)=\sigma^{\ord(c-d)}$ is a complete ultrametric space \cite{Berstel-Reutenauer_88}.
In this section two key lemmas are presented which describe how the distance between two series is
altered by operations involving the shuffle product defined in \eqref{shuffleproduct}. These results are employed in
Section~\ref{sec:transformation-group-allseriesdelta} to prove the existence of a group inverse
using ultrametric contractions on $(\allseries,\dist)$.

\begle
For any series $c_i,d_i\in\allseries$, $i=1,2$,
\begdi
\dist(c_1\shuffle d_1,c_2\shuffle d_2)\leq \max(\sigma^{\ord(c_1)}\dist(d_1,d_2),\sigma^{\ord(d_2)}\dist(c_1,c_2)).
\enddi
\endle
\begpr
First observe that
\begdi
\dist(c_i\shuffle d_1,c_i\shuffle d_2)=\sigma^{\ord(c_i\shuffle(d_1-d_2))}
=\sigma^{\ord(c_i)+\ord(d_1-d_2)}
=\sigma^{\ord(c_i)}\dist(d_1,d_2).
\enddi
In which case, it follows that
\begin{align*}
\dist(c_1\shuffle d_1,c_2\shuffle d_2)&\leq \max(\dist(c_1\shuffle d_1,c_1\shuffle d_2),\dist(c_1\shuffle d_2,c_2\shuffle d_2)) \\
&= \max(\sigma^{\ord(c_1)}\dist(d_1,d_2),\sigma^{\ord(d_2)}\dist(c_1,c_2)).
\end{align*}
\endpr

\begco \label{co:shuffle-nonexpansive}
For a fixed $c\in\allseries$, the mapping $d\mapsto c\shuffle d$ is an ultrametric contraction if $c$ is proper and an
isometry on $(\allseries,\dist)$ otherwise.
\endco

\begth {\rm{\cite{Gray-etal_Auto14}}} \label{th:shuffle-group}
The set of nonproper series $\allseriesnp\subset\allseries$ is a group under the shuffle
product. In particular, the shuffle inverse of any such series $c \in \allseriesnp$ is
$c^{\shuffle -1}=((c,\emptyset)(1-c^\prime))^{\shuffle -1}:=(c,\emptyset)^{-1}(c^{\prime})^{\shuffle\ast}$,
where $c^\prime:=1-c/(c,\emptyset)$ is proper and $(c^\prime)^{\shuffle\ast}:=\sum_{k\geq 0} (c^\prime)^{\shuffle k}$.
\endth

\begle \label{le:shuffle-inverse-isometry}
The shuffle inverse is an isometry for any $c,d\in\allseriesnp$ having identical
constant terms.
\endle
\begpr
For any $c,d\in\allseriesnp$ with $(c,\emptyset)=(d,\emptyset)$ observe that
\begdi
\ord(c^{\shuffle -1}-d^{\shuffle -1})=\ord\left(\sum_{k=1}^\infty (c^\prime)^{\shuffle k}-(d^\prime)^{\shuffle k}\right) \\
=\ord(c^\prime-d^\prime)=\ord(c-d),
\enddi
and hence the lemma is proved.
\endpr

\subsection{Hopf Algebras}
\label{subsec:HA}
Some definitions and basic facts concerning Hopf algebras are summarized in this section for later use.
The treatment is based largely on \cite{Figueroa-Gracia-Bondia_05,Hochschild_81,Loday-Ronco_10,Radford_12,Sweedler_69}.
In the following, a {\it{unital $\re$-algebra}} refers to a vector space $A$ over the base field $\mathbb{R}$ with associative product $m_A: A \otimes A \to A$ and unit map (element) $\mathsf{e}_A: \re \to A$ ($\mathsf{e}_A(1)=1_A$).
A {\it{counital coalgebra}} over $\re$ consists of a triple $(C,\Delta_C,\varepsilon_C)$. The coproduct $\Delta_C: C \to C \otimes C$ is coassociative, that is, $(\id_C \otimes \Delta_C)\circ \Delta_C=(\Delta_C \otimes \id_C)\circ \Delta_C$, where $\id_C$ is the identity map on $C$ and $\varepsilon_C: C \to \re$ denotes the counit (augmentation) map.
A coalgebra is cocommutative if $\tau_C \circ \Delta_C = \Delta_C$, where $\tau_C : C \otimes C \to C \otimes C$ is the flip map, $\tau_C(x \otimes y) = y \otimes x$. A {\it{bialgebra}} $B$ is both a unital algebra and a counital coalgebra together with compatibility relations, such as both the algebra product, $xy:=m_B(x\otimes y)$, and the unit map, $\mathsf{e}_B: \re \to B$, are coalgebra morphisms. This provides, for example, the identity $\Delta_B(xy)=\Delta_B(x)\Delta_B(y)$. The unit of $B$ is denoted by $\mathbf{1}_B = \mathsf{e}_B(1)$. A bialgebra is called {\it{graded}} if there are $\re$-vector subspaces $B_n$, $n \in \mathbb{N}$ such that $B= \bigoplus_{n \geq 0} B_n$ with $m_B(B_k \otimes B_l) \subseteq B_{k+l}$ and $\Delta_B B_n \subseteq \bigoplus_{k+l=n} B_k\otimes B_l.$ Elements $x \in B_n$ are given a degree
$\deg(x)=n$. Moreover, $B$ is called {\it{connected}} if the $\re$-vector subspace $B_0 = \re\mathbf{1}_B$.
Define $B_+=\bigoplus_{n > 0} B_n$. In a connected graded bialgebra, the coproduct for any $x \in B_n$ is of the form
\begdi
	\Delta_B(x) = x \otimes \mathbf{1}_B + \mathbf{1}_B \otimes x + \Delta_B'(x) \in \bigoplus_{k+l=n} B_k \otimes B_l,
\enddi
where $\Delta_B'(x) := \Delta_B(x) -  x \otimes \mathbf{1}_B - \mathbf{1}_B \otimes x \in \bigoplus_{k+l=n \atop k,l>0} B_k \otimes B_l$ is the {\em reduced} coproduct. Note that in the following the use of subscripts on the structure maps is limited to the cases where there is potential for notational confusion.

Suppose $B$ is a bialgebra and $A$ is an $\re$-algebra with product $m_A: A \otimes A \to A$ and unit map $\mathsf{e}_A$, e.g., $A=\re$ or $A=B$. The vector space  $L(B, A)$ of linear maps from the bialgebra $B$ to $A$ together with the convolution product
\begin{equation}
\label{convol}
	\Phi \star \Psi := m_{A} \circ (\Phi \otimes \Psi) \circ \Delta : B \to A,
\end{equation}
where $\Phi,\Psi \in L(B,A)$, is an associative algebra with unit $\iota := \mathsf{e}_{A} \circ \varepsilon$.
A {\it{Hopf algebra}} $H$ is a bialgebra together with a particular $\re$-linear map called an {\it{antipode}} $S: H \to H$
which satisfies the Hopf algebra axioms
\cite{Radford_12,Sweedler_69} and has the property of being an antihomomorphism for both the algebra and the coalgebra structures, i.e., $S(xy)=S(y)S(x)$ and $\Delta \circ S = (S \otimes S)\circ \tau \circ \Delta$. When $A=H$, the necessarily unique antipode $S \in L(H,H)$ is the inverse of the identity map $\id : H \to H$ with respect to the convolution product, that is,
\begin{equation}
\label{antipode}
    S  \star\id = \id \star S := m \circ (S \otimes \id) \circ \Delta = \mathsf{e} \circ \varepsilon.
\end{equation}
{\it{Group-like}} and {\it{primitive}} elements in a Hopf algebra satisfy $\Delta(x)=x \otimes x$ and $\Delta(y) = y \otimes \mathbf{1} + \mathbf{1} \otimes y$, respectively. Note that if $x_1,x_2$ are primitive in $H$, then $[x_1,x_2]:=m_H(x_1 \otimes x_2) - m_H( x_2 \otimes x_1)$ is primitive as well. In particular, the set of primitive elements $P(H) \subset H$ is a Lie algebra.
Another important observation is that a connected graded bialgebra $H=\bigoplus_{n \ge 0} H_n$ is {\em always} a connected graded Hopf algebra.

Let $H = \bigoplus_{n \geq 0} H_n$ be a connected graded commutative bialgebra, and suppose $A$ is a commutative unital algebra. Any $\Phi \in L(H, A)$ is called a {\it{character}} if $\Phi(\mathbf{1})=1_A$ and
\begin{equation}
\label{character}
	\Phi(xy) = \Phi(x)\Phi(y)
\end{equation}
for all $x,y\in H$. The set of characters is denoted by $G_{ A} \subset L(H, A)$ and forms a group with respect to the convolution product \eqref{convol}. The neutral element $\iota:=\mathsf{e}_{ A}\circ \varepsilon$ in $G_{ A}$ is given by $\iota(\mathbf{1})=1_{ A}$ and $\iota(x) = 0$ for elements $x$  in the {\it{augmentation ideal}} $\mathsf{Ker}(\varepsilon)=H_+=\bigoplus_{n > 0} H_n$.
The inverse of $\Phi \in G_{ A}$ is given by composition with the antipode
\begeq
\label{eq:convinv-via-antipode}
	\Phi^{\star -1} = \Phi \circ S,
\endeq
which follows from \eqref{antipode} and the fact that $\Phi$ is an algebra morphism. Specifically,
\begin{align*}
	\Phi \star (\Phi \circ S) &= \Phi \circ (\id \star S) = \Phi\circ (S \star \id)= (\Phi\circ S) \star \Phi
					  = \Phi\circ \mathsf{e} \circ \varepsilon=\mathsf{e} \circ \varepsilon.
\end{align*}
An {\it{infinitesimal character}} with values in $A$ is a linear map $\xi \in L(H, A)$ such that for $x, y \in  H^+$, $\xi(m_H(x \otimes y)) = 0$, which implies $\xi(\mathbf{1}) =0$. The linear space of infinitesimal characters $g_{A} \subset L(H, A)$ forms a Lie algebra with respect to the Lie bracket defined in terms of the convolution product \eqref{convol}. For $\alpha \in g_{A}$ and any $x \in H_n$ the exponential $\exp^\star(\alpha)(x) := \sum_{j\ge 0} \frac{1}{j!}\alpha^{\star j}(x)$ is a finite sum terminating at $j=n$.
\begin{proposition}\label{prop:exp}
$\exp^\star$ restricts to a bijection from $g_A$ onto $G_A$.
\end{proposition} \\
(For details see references \cite{EbraManchon06, Figueroa-Gracia-Bondia_05}). The inverse of $\exp^\star$ is given by $\log^\star(\iota+\gamma)(x)=\sum_{k\ge 1}{(-1)^{k-1}\over k}\gamma^{\star k}(x)$, where again the sum terminates at $k=n$ for any $x \in H_n$.

Finally, a brief description of the relationship between the group $G_A \subset L(H, A)$ and the Hopf algebra $H$ is given. The reader is referred to the paper of Manchon and Frabetti \cite{Frabetti-Manchon_15} for details and additional references. By definition, elements in $G_A$ map all of $H$ into the commutative unital algebra $A$. However, any element $x \in H$ can also be seen as an $A$-valued function on $G_A$. Indeed, let $\Phi \in G_A$, then $x(\Phi):=\Phi(x) \in A$ and the usual pointwise product of functions
$
	(xy)(\Phi)=x(\Phi)y(\Phi)
$
follows from \eqref{character} since $\Phi \in G_A$. The definition of the convolution product \eqref{convol} in terms of the coproduct of $H$ implies a natural coproduct on functions $x \in H$, that is, $\Delta(x)(\Phi,\Psi) :=(\Phi \star \Psi)(x) \in A$. Similarly, the inverse of $G_A$ as well as its unit  correspond naturally to the antipode and counit map on $H$, respectively. This {\em reversed} perspective on the relationship between $H$ and its group of characters $G_A$ allows one to interpret $H$ as the (Hopf) algebra of  {\it{coordinate functions}} of the group $G_A$.
More precisely, $H$ contains the {\em representative functions} over $G_A$. The reader is directed to Cartier's work \cite{Cartier_07} for a comprehensive review of this topic. In this paper, the starting point will be a particular group whose product, unit and inverse are used to identify its Hopf algebra of coordinate functions.

\section{Output Affine Feedback Transformation Group}
\label{sec:transformation-group-allseriesdelta}

In general, if $G$ is a group and $S$ a given set, then $G$ is said to act as a {\em transformation
group} on the right of $S$ if there exists a mapping
$\phi:S\times G\rightarrow S:(h,g)\mapsto hg$
such that:
\begin{description}
\item[i.] $h1=h$, where $1$ is the identity element of $G$;
\item[ii.] $h(g_1g_2)=(h g_1)g_2$ for all $g_1,g_2\in G$.
\end{description}
The action is said to be {\em free} if $hg=h$ implies that $g=1$.
It will be evident in what follows that it is most convenient to redefine $\allseriesdelta:=\allseries \times \allseries$ and to identify the output affine feedback transformation group with those
elements of $\allseriesdelta$ which are invertible in some sense. Each tuple $c_\delta:=(c_L,c_R) \in \allseriesdeltaLC:=\allseriesLC \times \allseriesLC$ gives rise to an operator of the form $F_{c_{\delta}}[u]:=F_{c_L}[u]u+F_{c_R}[u]$,
and by design the group product $c_\delta \circ d_\delta$ for elements in $\allseriesdeltaLC$ will satisfy the identity $F_{c_\delta} \circ F_{d_\delta}=F_{c_\delta\circ d_\delta}$
when $c_\delta\circ d_\delta\in\allseriesdeltaLC$.
In addition, a second product $c \modcomp d$ is defined which describes the action of the group $\allseriesdelta$ on $\allseries$ so that $F_c \circ F_{d_\delta} = F_{c\modcomp d_\delta}$ when the series involved are all locally convergent.
This product, referred to as the {\em mixed composition product},
plays a central role in the paper, so it will be the starting point for the section.
Once these concepts are in place and the transformation group is established, the affine
extension of the output feedback product in Theorem~\ref{th:feedback-product-formula}
will be given.
Henceforth, the focus is on the SISO case.

For a fixed $d_{\delta}=(d_L,d_R)\in\allseriesdelta$, the mixed composition product is defined in
terms of an algebra morphism $\phi_d$ which takes words in $X^\ast$ to a composition of simple vector space endomorphisms involving only the shuffle product and
catenation with a single letter $x_0$ or $x_1$. This sets up an inductive definition which provides some
insight into the core combinatorial structures at play. The reader is referred in particular to \cite{Foissy_CiA15,Foissy_EJM15} for a
broader theoretical view behind this approach.

\begde \label{de:redefined-allseriesdelta}
Let $d_\delta=(d_L,d_R)\in\allseriesdelta$.
Define the \bfem{mixed composition product} mapping $\allseries\times\allseriesdelta$ into
$\allseries$ as
\begin{align}
	c\modcomp d_\delta &=\phi_d(c)(1)
					=\sum_{\eta\in X^\ast} (c,\eta)\, \phi_d(\eta)(1), \label{def:mixedCompProd}
\end{align}
where $\phi_d$ is the continuous (in the ultrametric sense) algebra homomorphism from $\allseries$ to $\Endallseries$ uniquely specified by
$\phi_d(x_i\eta)=\phi_d(x_i)\circ \phi_d(\eta)$ with
\begeq \label{eq:phi-definition}
\phi_d(x_0)(e)=x_0e,\;\;
\phi_d(x_1)(e)=x_1(d_L\shuffle e)+x_0(d_R\shuffle e)
\endeq
for any $e\in\allseries$, and where $\phi_d(\emptyset)$ denotes the identity map on $\allseries$.
\endde

\begex
Suppose $c=2x_0+x_1x_0$ and $d_\delta=(d_L,d_R)=(3x_1,-x_1^2)$. Then
\begdi
c\modcomp d_{\delta}=2\phi_d(x_0)(1)+\phi_d(x_1)\circ\phi_d(x_0)(1),
\enddi
where $\phi_d(x_0)(e)=x_0e$ and
$\phi_d(x_1)(e)=x_1(3x_1\shuffle e)+x_0(-x_1^2\shuffle e)$.
Therefore,
$c\modcomp d_{\delta}=2x_0+3x_1^2x_0+3x_1x_0x_1-x_0x_1^2x_0-x_0x_1x_0x_1-x_0^2x_1^2$.
\endex

The modified composition product referred to in Section~\ref{sec:preliminaries}
and used in earlier work on output feedback \cite{Gray-etal_Auto14,Gray-etal_SCL14,Gray-Li_05,Li_04,Thitsa-Gray_12}
corresponds here to the special case where $d_L=1$.
Some fundamental properties of the mixed composition product are given next.

\begle \label{le:mixed-product-properties}
The mixed composition product \eqref{def:mixedCompProd}
\begin{description}
\item[(1)] is left $\re$-linear;
\item[(2)] satisfies $c\modcomp (1,0)=c$;
\item[(3)] satisfies $c\modcomp d_\delta=k\in\re$ for any fixed $d_{\delta}$ if and only if $c=k$;
\item[(4)] satisfies $(x_0c)\modcomp d_{\delta}=x_0(c\modcomp d_{\delta})$ and $(x_1c)\modcomp d_{\delta}=x_1(d_L\shuffle (c\modcomp d_{\delta}))+x_0(d_R\shuffle (c\modcomp d_{\delta}))$;
\item[(5)] distributes to the left over the shuffle product.
\end{description}
\endle
\begpr
$\;$\\
\noindent(1) This fact follows directly from the definition of the mixed composition product.
\newline
\noindent(2) The claim is immediate since $\phi_{(1,0)}(\eta)(1)=\eta$.
\newline
\noindent(3) The only non trivial assertion is that $c\modcomp d_{\delta}=k$ implies $c=k$.
This claim is best handled in a Hopf algebra setting. So the proof is deferred until Section~\ref{sec:hopf-algebra}.
\newline
\noindent(4) Observe:
\begin{align*}
(x_0c)\modcomp d_\delta&=\phi_d(x_0c)(1)=\phi_d(x_0)\circ \phi_d(c)(1)=x_0(c\modcomp d_{\delta}) \\
(x_1c)\modcomp d_\delta&=\phi_d(x_1c)(1)=\phi_d(x_1)\circ \phi_d(c)(1)
=x_1(d_L\shuffle(c\modcomp d_{\delta}))+x_0(d_R\shuffle(c\modcomp d_{\delta})).
\end{align*}
\noindent(5) For any $e_\delta\in\allseriesdelta$, one can define a shuffle product on
$\Endallseries$ via
\begdi
\phi_e(x_i\eta)\shuffle\phi_e(x_j\xi)=
\phi_e(x_i)\circ [\phi_e(\eta)\shuffle\phi_e(x_j\xi)]+
\phi_e(x_j)\circ [\phi_e(x_i\eta)\shuffle \phi_e(\xi)].
\enddi
In which case, $\phi_e$ acts as an algebra map between the shuffle algebra on $\allseries$ and
the shuffle algebra on $\Endallseries$. That is,
$
\phi_e(c\shuffle d)=\phi_e(c)\shuffle \phi_e(d).
$
Hence, $(c\shuffle d)\modcomp e_\delta=\phi_e(c\shuffle d)(1)=\phi_e(c)(1)\shuffle \phi_e(d)(1)=(c\modcomp e_\delta)\shuffle(d\modcomp e_\delta)$.
\endpr

It is easily checked that
$
\dist(c_\delta,d_\delta):=\max(\dist(c_L,d_L),\dist(c_R,d_R))
$
is an ultrametric on $\allseriesdelta$, and $\allseriesdelta$ is complete.\footnote{Using $\dist$
for both the ultrametric on $\allseries$ and $\allseriesdelta$
should cause minimal confusion since their arguments are distinct.}
The following lemma states that the
mixed composition product reduces ultrametric distance in a certain sense.

\begle \label{le:modcomp-contraction-on-group}
For any $c\in\allseries$ and $d_{\delta,1},d_{\delta,2}\in\allseriesdelta$ it follows that
\begdi
\dist(c\modcomp d_{\delta,1},c\modcomp d_{\delta,2})\leq \sigma^{\ord(c^\prime)}\dist(d_{\delta,1},d_{\delta,2}),
\enddi
where $c=(c,\emptyset)\emptyset+c^\prime$ with $c^\prime$ proper.
\endle

\begpr
For a fixed $d_R$, consider the map $d_L\mapsto c\modcomp (d_L,d_R)$. Likewise, for a
fixed $d_L$ there is a companion map
$d_R\mapsto c\modcomp (d_L,d_R)$. It is first shown that
on the ultrametric space $(\allseries,\dist)$:
\begin{subequations} \label{eq:modcomp-contraction-dx-fixed}
\begin{align}
\dist(c\modcomp (d_{L,1},d_R),c\modcomp (d_{L,2},d_R))&\leq \sigma^{\ord(c^\prime)}\dist(d_{L,1},d_{L,2}) \label{eq:modcomp-contraction-dR-fixed} \\
\dist(c\modcomp (d_{L},d_{R,1}),c\modcomp (d_{L},d_{R,2}))&\leq \sigma^{\ord(c^\prime)}\dist(d_{R,1},d_{R,2}). \label{eq:modcomp-contraction-dL-fixed}
\end{align}
\end{subequations}
Inequality \rref{eq:modcomp-contraction-dR-fixed} is clearly true when $d_{L,1}=d_{L,2}$. When $d_{L,1}\neq d_{L,2}$,
it is first necessary to show that for any $\eta\in X^\ast$
\begeq \label{eq:ord-phid1-minus-phid2}
\ord(\phi_{d_1}(\eta)(1)-\phi_{d_2}(\eta)(1))\geq \abs{\eta}+\ord(d_{L,1}-d_{L,2}).
\endeq
The proof is by induction on the length of $\eta$.
The claim is trivial when $\eta$ is empty or a single letter. Assume the inequality holds for words up to length
$k\geq 0$. For any $x_0\eta$ with $\eta\in X^k$, inequality \rref{eq:ord-phid1-minus-phid2} follows directly from
the induction hypothesis.
The case for $x_1\eta$ is handled as follows:
\begin{align*}
\lefteqn{\ord(\phi_{d_1}(x_1\eta)(1)-\phi_{d_2}(x_1\eta)(1))}\hspace*{0.2in} \\
&=\ord(x_1[d_{L,1}\shuffle \phi_{d_1}(\eta)(1)-d_{L,2}\shuffle \phi_{d_2}(\eta)(1)]+
x_0[d_R\shuffle (\phi_{d_1}(\eta)(1)- \phi_{d_2}(\eta)(1))]) \\
&=\ord(x_1[(d_{L,1}-d_{L,2})\shuffle \phi_{d_1}(\eta)(1)+
d_{L,2}\shuffle (\phi_{d_1}(\eta)(1)-\phi_{d_2}(\eta)(1))]+ \\
&\hspace*{0.18in}x_0[d_R\shuffle (\phi_{d_1}(\eta)(1)- \phi_{d_2}(\eta)(1))])\\
&\geq 1 + \min(\ord([d_{L,1}-d_{L,2}]\shuffle \phi_{d_1}(\eta)(1)),
\ord(d_{L,2}\shuffle [\phi_{d_1}(\eta)(1)-\phi_{d_2}(\eta)(1)]), \\
&\hspace*{0.18in}\ord(d_R\shuffle [\phi_{d_1}(\eta)(1)- \phi_{d_2}(\eta)(1)]))\\
&\geq 1 + \min(\ord(d_{L,1}-d_{L,2})+\abs{\eta},\ord(\phi_{d_1}(\eta)(1)-
\phi_{d_2}(\eta)(1)))\\
&= \abs{\eta}+1 + \ord(d_{L,1}-d_{L,2}).
\end{align*}
In which case, \rref{eq:ord-phid1-minus-phid2} holds for any $\eta\in X^\ast$.
The inequality \rref{eq:modcomp-contraction-dR-fixed} is now derived. Observe
\begin{align*}
\lefteqn{\dist(c\modcomp(d_{L,1},d_R),c\modcomp  (d_{L,2},d_R))} \hspace*{0.5in}\\
&=\dist(c^\prime \modcomp (d_{L,1},d_R),c^\prime\modcomp (d_{L,2},d_R))
=\sigma^{\ord(\sum_{\eta} (c^{\prime},\eta)[\phi_{d_1}(\eta)(1)-\phi_{d_2}(\eta)(1)])} \\
&\leq\sigma^{\min_{\eta\in\supp(c^\prime)}\ord(\phi_{d_1}(\eta)(1)-\phi_{d_2}(\eta)(1))}
\leq \sigma^{\min_{\eta\in\supp(c^\prime)}\abs{\eta}+\ord(d_{L,1}-d_{L,2})} \\
&=\sigma^{\ord(c^\prime)}\dist(d_{L,1},d_{L,2}).
\end{align*}
The proof for
\rref{eq:modcomp-contraction-dL-fixed} is completely analogous.
The final step of the proof is to employ the ultrametric property {\rm iv} in
conjunction with \rref{eq:modcomp-contraction-dx-fixed}.
Observe
\begin{align*}
\lefteqn{\dist(c\modcomp d_{\delta,1},c\modcomp d_{\delta,2})} \hspace*{0.2in}\\
&=\dist(c\modcomp (d_{L,1},d_{R,1}),c\modcomp (d_{L,2},d_{R,2})) \\
&\leq \max(\dist(c\modcomp (d_{L,1},d_{R,1}),c\modcomp (d_{L,2},d_{R,1})),
\dist(c\modcomp (d_{L,2},d_{R,1}),c\modcomp (d_{L,2},d_{R,2}))) \\
&\leq \sigma^{\ord(c^\prime)} \max(\dist(d_{L,1},d_{L,2}),\dist(d_{R,1},d_{R,2}))
=\sigma^{\ord(c^\prime)} \dist(d_{\delta,1},d_{\delta,2}).
\end{align*}
\endpr

The group product on $\allseriesdelta$ is now defined in terms of the mixed composition product.
Its basic properties are given in the subsequent lemma.

\begde \label{de:extended-feedback-group-product}
The \bfem{composition product} on $\allseriesdelta$ is defined as
\vspace*{-0.05in}
\begdi
c_\delta\circ d_\delta = ((c_L\modcomp d_\delta)\shuffle d_L,(c_L\modcomp d_\delta)\shuffle d_R+c_R\modcomp d_\delta).
\enddi
\endde

\begle \label{le:group-product-identities}
The composition product on $\allseriesdelta$
\begin{description}
\item[(1)] is left $\re$-linear;
\item[(2)] satisfies $(c\modcomp d_\delta)\modcomp e_\delta=c\modcomp (d_\delta\circ e_\delta)$ (mixed associativity);
\item[(3)] is associative.
\end{description}
\endle
\begpr
$\;$\\
\noindent(1) This property follows from the left linearity of the mixed composition product.
\newline
\noindent(2) In light of the first item, it is sufficient to prove the claim only for $c=\eta\in X^k$, $k\geq 0$.
The cases $k=0$ and $k=1$ are trivial.  Assume the claim holds up to some fixed $k\geq 0$.
Then via Lemma~\ref{le:mixed-product-properties} (4) and the induction hypothesis it follows that
\begdi
((x_0\eta)\modcomp d_\delta)\modcomp e_\delta=(x_0(\eta\modcomp d_\delta))\modcomp e_\delta
=x_0((\eta\modcomp d_\delta)\modcomp e_\delta)
=x_0(\eta\modcomp (d_\delta\circ e_\delta))
=(x_0\eta)\modcomp (d_\delta\circ e_\delta).
\enddi
In a similar fashion, apply the properties in Lemma~\ref{le:mixed-product-properties} (1), (4), and (5) to get
\begin{align*}
\lefteqn{((x_1\eta)\modcomp d_\delta)\modcomp e_\delta} \hspace*{0.0in}\\
&=[x_1(d_L\shuffle (\eta\modcomp d_{\delta}))+x_0(d_R\shuffle (\eta\modcomp d_{\delta}))]\modcomp e_\delta \\
&=[x_1(d_L\shuffle (\eta\modcomp d_{\delta}))]\modcomp e_\delta+[x_0(d_R\shuffle (\eta\modcomp d_{\delta}))]\modcomp e_\delta \\
&=x_1[e_L\shuffle ((d_L\shuffle (\eta\modcomp d_{\delta}))\modcomp e_\delta)]+
x_0[e_R\shuffle ((d_L\shuffle (\eta\modcomp d_{\delta}))\modcomp e_\delta)]+ \\
&\hspace*{0.18in}x_0[(d_R\shuffle (\eta\modcomp d_{\delta}))\modcomp e_\delta] \\
&=x_1[\underbrace{e_L\shuffle (d_L\modcomp e_\delta)}_{(d_\delta\circ e_\delta)_L}\shuffle ((\eta\modcomp d_{\delta})\modcomp e_\delta)]+
x_0[\underbrace{((d_L\modcomp e_\delta)\shuffle e_R+d_R\modcomp e_\delta)}_{(d_\delta \circ e_\delta)_R}\shuffle ((\eta\modcomp d_{\delta})\modcomp e_\delta)].
\end{align*}
Now employ the induction hypothesis so that
\begin{align*}
((x_1\eta)\modcomp d_\delta)\modcomp e_\delta
&=x_1[(d_\delta\circ e_\delta)_L\shuffle(\eta\modcomp(d_\delta\circ e_\delta))]+
x_0[(d_\delta\circ e_\delta)_R\shuffle (\eta\modcomp(d_\delta\circ e_\delta))] \\
&=(x_1\eta)\modcomp (d_\delta\circ e_\delta).
\end{align*}
Therefore, the claim holds for all $\eta\in X^\ast$, and the identity is proved.
\newline
\noindent(3) First apply Definition~\ref{de:extended-feedback-group-product} twice, Lemma~\ref{le:mixed-product-properties} (1) and (5) to get
\begin{align*}
(c_\delta\circ d_\delta)\circ e_\delta
&=((c_L\modcomp d_\delta)\shuffle d_L,(c_L\modcomp d_\delta)\shuffle d_R+c_R\modcomp d_\delta)\circ e_\delta \\
&=([((c_L\modcomp d_\delta)\shuffle d_L)\modcomp e_\delta]\shuffle e_L,[((c_L\modcomp d_\delta)\shuffle d_L)\modcomp e_\delta]\shuffle e_R+ \\
&\hspace*{0.18in}[(c_L\modcomp d_\delta)\shuffle d_R+c_R\modcomp d_\delta]\modcomp e_\delta) \\
&=([(c_L\modcomp d_\delta)\modcomp e_\delta]\shuffle [d_L\modcomp e_\delta]\shuffle e_L,
[(c_L\modcomp d_\delta)\modcomp e_\delta]\shuffle \\
&\hspace*{0.18in}[d_L\modcomp e_\delta]\shuffle e_R+
((c_L\modcomp d_\delta)\modcomp e_\delta) \shuffle (d_R\modcomp e_\delta)+(c_R\modcomp d_\delta)\modcomp e_\delta).
\end{align*}
Now apply the mixed associativity property from the previous item and then recombine terms according to Definition~\ref{de:extended-feedback-group-product}
so that
\begin{align*}
(c_\delta\circ d_\delta)\circ e_\delta
&=([c_L\modcomp (d_\delta\circ e_\delta)]\shuffle \underbrace{[d_L\modcomp e_\delta]\shuffle e_L}_{(d_\delta\circ e_\delta)_L},
[c_L\modcomp (d_\delta\circ e_\delta)]\shuffle
[d_L\modcomp e_\delta]\shuffle e_R+ \\
&\hspace*{0.18in}(c_L\modcomp (d_\delta\circ e_\delta)) \shuffle (d_R\modcomp e_\delta)+
c_R\modcomp (d_\delta\circ e_\delta)) \\
&=([c_L\modcomp (d_\delta\circ e_\delta)]\shuffle (d_\delta\circ e_\delta)_L,
[c_L\modcomp (d_\delta\circ e_\delta)]\shuffle \\
&\hspace*{0.18in}
[\underbrace{(d_L\modcomp e_\delta)\shuffle e_R+d_R\modcomp e_\delta}_{(d_\delta\circ e_\delta)_R}]
+c_R\modcomp (d_\delta\circ e_\delta)) \\
&= ((c_\delta\circ (d_\delta\circ e_\delta))_L,(c_\delta\circ (d_\delta\circ e_\delta))_R)
=c_\delta\circ(d_\delta\circ e_\delta),
\end{align*}
and the lemma is proved.
\endpr

The Fliess operator interpretations of the two formal power series products described above are given in the next theorem.

\begth
For any $c\in\allseriesLC$ and $c_\delta,d_\delta\in\allseriesdeltaLC$ such that
$c\modcomp d_{\delta}\in\allseriesLC$ and $c_\delta\circ d_\delta\in \allseriesdeltaLC$,
the following identities hold:
\begdes
\item[(1)] $F_c\circ F_{d_{\delta}}=F_{c\modcomp d_{\delta}}$
\item[(2)] $F_{c_{\delta}}\circ F_{d_{\delta}}= F_{c_\delta\circ d_\delta}$.
\enddes
\endth

\begpr
$\;$ \\
\noindent (1) It is sufficient to prove the claim for $c=\eta\in X^\ast$. This is done by induction on the length of $\eta$.
The case for the empty word is trivial. Assume the identity holds for words $\eta\in X^k$ up to some fixed length $k\geq 0$. Then
\begin{align*}
E_{x_0\eta}\circ F_{d_\delta}[u](t)&=\int_{t_0}^t E_{\eta}[F_{d_\delta}[u]](\tau,t_0)\,d\tau
=\int_{t_0}^t F_{\eta\modcomp d_\delta}[u](\tau)\,d\tau
=F_{x_0(\eta\modcomp d_{\delta})}[u](t) \\
&=F_{(x_0\eta)\modcomp d_{\delta}}[u](t).
\end{align*}
Similarly,
\begin{align*}
E_{x_1\eta}\circ F_{d_\delta}[u](t)
&=\int_{t_0}^t F_{\eta\modcomp d_\delta}[u](\tau)[F_{d_L}[u](\tau)u(\tau)+F_{d_R}[u](\tau)]\,d\tau \\
&=\int_{t_0}^t F_{d_L\shuffle (\eta\modcomp d_{\delta})}[u](\tau)u(\tau)+F_{(\eta\modcomp d_{\delta})\shuffle d_R}[u](\tau)\,d\tau \\
&=F_{x_1(d_L\shuffle (\eta\modcomp d_{\delta}))+x_0((\eta\modcomp d_{\delta})\shuffle d_R)}[u](t)
=F_{(x_1\eta)\circ d_\delta}[u](t).
\end{align*}
Hence, the claim holds for all $\eta\in X^\ast$.

\noindent (2) Observe
\begin{align*}
F_{c_\delta}\circ F_{d_{\delta}}[u]&=(F_{c_L}[u]u+F_{c_R}[u])\circ (F_{d_L}[u]u+F_{d_R}[u]) \\
&=F_{c_L}[F_{d_{\delta}}[u]]F_{d_L}[u]u+F_{c_L}[F_{d_{\delta}}[u]]F_{d_R}[u]+
F_{c_R}[F_{d_\delta}[u]] \\
&=F_{(c_L\modcomp d_{\delta})\shuffle d_L}[u]u+F_{(c_L\modcomp d_{\delta})\shuffle d_R}[u]+F_{c_R\modcomp d_\delta}[u] \\
&=F_{(c_\delta\circ d_\delta)_L}[u]u+F_{(c_\delta\circ d_\delta)_R}[u]
=F_{c_\delta\circ d_\delta}[u].
\end{align*}

\endpr

Let $\allseriesdeltanp$ denote the subset of $\allseriesdelta$ whose left series are nonproper.

\begth
The set $(\allseriesdeltanp,\circ,(1,0))$ is a group.
\endth

\begpr
From Lemma~\ref{le:mixed-product-properties} (2), it follows directly that $c_\delta\circ (1,0)=c_\delta$.
Using the identities $1\modcomp c_\delta=1$ and $0\modcomp c_\delta=0$ from Lemma~\ref{le:mixed-product-properties} (3), it is straightforward to see that
$(1,0)\circ c_\delta=c_\delta$.  Associativity of the group product was established in
Lemma~\ref{le:group-product-identities} (3). So the only open issue
is the existence of inverses. Suppose $c_\delta$ is fixed, and one seeks a right inverse $c^{-1}_\delta:=(c_L^{\circ -1},c_R^{\circ -1})$, that is,
$c_\delta\circ c^{-1}_\delta=(1,0)$.\footnote{The notation $c^{\circ -1}$, $c\in\allseries$, is used in this situation instead of $c^{-1}$
to distinguish it from $c^{\shuffle -1}$.}
Then it follows directly from Theorem~\ref{th:shuffle-group} and
Definition~\ref{de:extended-feedback-group-product} that
\begin{subequations}
\label{eq:right-inverse-generalized-feedback-inverse}
\begin{align}
c_L^{\circ -1}&=(c_L\modcomp (c_L^{\circ -1},c_R^{\circ -1}))^{\shuffle -1} \label{eq:cL-right-inverse} \\
c_R^{\circ -1}&=-c_L^{\circ -1}\shuffle (c_R\modcomp (c_L^{\circ -1},c_R^{\circ -1})). \label{eq:cR-right-inverse}
\end{align}
\end{subequations}
It is first shown that the mapping
\begdi
S_R\!:\!(e_L,e_R)\mapsto\! ((c_L\modcomp \!(e_L,e_R))^{\shuffle -1},-e_L\shuffle \!(c_R\modcomp \!(e_L,e_R)))
\enddi
is an ultrametric contraction on $\allseriesdelta$, and therefore has a unique fixed point, which by design is a right inverse, $c^{-1}_\delta$.
Note that for any $e_\delta$ it follows that $(S_R(e_L,e_R)_L,\emptyset)=(c_L,\emptyset)^{-1}\neq 0$. Thus, the fixed point will
always be in the group.
Then it is shown that this same series is also
a left inverse, that is, $c_\delta^{-1}\circ c_\delta=(1,0)$, or equivalently,
\begin{subequations}
\label{eq:left-inverse-generalized-feedback-inverse}
\begin{align}
c_L&=(c_L^{\circ -1}\modcomp (c_L,c_R))^{\shuffle -1}  \label{eq:cL-left-inverse}\\
c_R&=-c_L\shuffle(c_R^{\circ -1}\modcomp (c_L,c_R)). \label{eq:cR-left-inverse}
\end{align}
\end{subequations}
To establish the first claim, observe via Corollary~\ref{co:shuffle-nonexpansive} and
Lemma~\ref{le:shuffle-inverse-isometry} that for any $e_\delta,\bar{e}_\delta\in\allseriesdelta$
\begin{align*}
\dist(S_R(e_\delta),S_R(\bar{e}_\delta))
&=\max(
\dist((c_L\modcomp \!(e_L,e_R))^{\shuffle -1},(c_L\modcomp \!(\bar{e}_L,\bar{e}_R))^{\shuffle -1}), \\
&\hspace*{0.18in}\dist(-e_L\shuffle (c_R\modcomp (e_L,e_R)),-\bar{e}_L\shuffle (c_R\modcomp (\bar{e}_L,\bar{e}_R)))) \\
&\leq\max(
\dist(c_L\modcomp (e_L,e_R),c_L\modcomp (\bar{e}_L,\bar{e}_R)), \\
&\hspace*{0.18in}
\dist(c_R\modcomp (e_L,e_R),c_R\modcomp (\bar{e}_L,\bar{e}_R))).
\end{align*}
In which case, from Lemma~\ref{le:modcomp-contraction-on-group} it follows that
\begin{align*}
\dist(S_R(e_\delta),S_R(\bar{e}_\delta))
&\leq \max(\sigma^{\ord(c_L^\prime)}\dist((e_L,e_R),(\bar{e}_L,\bar{e}_R)), \\
&\hspace*{0.18in}
\sigma^{\ord(c_R^\prime)}\dist((e_L,e_R),(\bar{e}_L,\bar{e}_R)))\\
&\leq \sigma\, \dist(e_\delta,\bar{e}_\delta).
\end{align*}
To address the second claim, suppose $c_\delta^{-1}$ satisfies \rref{eq:cL-right-inverse}. In which case,
\begin{align*}
(c_L\modcomp c_{\delta}^{-1})\shuffle c_L^{\circ -1}&=1 \\
(c_L\modcomp c_{\delta}^{-1})\shuffle (c_L^{\circ -1}\modcomp (c_\delta\circ c_\delta^{-1}))&=1.
\end{align*}
Using Lemma~\ref{le:group-product-identities} (2) and Lemma~\ref{le:mixed-product-properties} (5) gives
\begin{align*}
(c_L\modcomp c_{\delta}^{-1})\shuffle ((c_L^{\circ -1}\modcomp c_\delta)\modcomp c_\delta^{-1}))&=1 \\
(c_L\shuffle (c_L^{\circ -1}\modcomp c_\delta))\modcomp c_\delta^{-1}&=1.
\end{align*}
Applying Lemma~\ref{le:mixed-product-properties} (3) then yields
$c_L\shuffle (c_L^{\circ -1}\modcomp c_\delta)=1$, and thus,
$c_L=(c_L^{\circ -1}\modcomp c_\delta)^{\shuffle -1}$, which is \rref{eq:cL-left-inverse}.
If, in addition, $c_\delta^{-1}$ also satisfies \rref{eq:cR-right-inverse}, then substituting \rref{eq:cL-right-inverse}
into \rref{eq:cR-right-inverse} gives $c_R^{\circ -1}=-(c_L\modcomp c_\delta^{-1})^{\shuffle -1}\shuffle (c_R\modcomp c_\delta^{-1})$.
Therefore, in a similar fashion,
\begin{align*}
-(c_L\modcomp c_\delta^{-1})\shuffle c_R^{\circ -1}&=c_R\modcomp c_\delta^{-1} \\
(-c_L\shuffle (c_R^{\circ -1}\modcomp c_\delta))\modcomp c_\delta^{-1}&=c_R\modcomp c_\delta^{-1} \\
(c_R+c_L\shuffle (c_R^{\circ -1}\modcomp c_\delta))\modcomp c_\delta^{-1}&=0.
\end{align*}
Once again applying Lemma~\ref{le:mixed-product-properties} (3) gives
$c_R+c_L\shuffle (c_R^{\circ -1}\modcomp c_\delta)=0$,
which is equivalent to \rref{eq:cR-left-inverse}.
\endpr

The following theorem is the first main result of this section. The issue of
freeness will be addressed in Section~\ref{sec:group-invariants}.

\begth
The group $(\allseriesdeltanp,\circ,(1,0))$ acts as a right transformation group on
$\allseries$ via the action $c\modcomp d_{\delta}$.
\endth
\begpr
See Lemma~\ref{le:mixed-product-properties} (2) and Lemma~\ref{le:group-product-identities} (2).
\endpr

The second main result of the section is a formula for the group element $e_\delta$ which
describes the output affine feedback connection shown in Figure~\ref{fig:affine-output-feedback-connection}.
This defines the {\em output affine feedback product}.
\begth
Let $c\in\allseriesLC$ and $d_{\delta}=(d_L,d_R)\in\allseriesdeltaLC$. Then the generating series for the output
affine feedback connected system is
$
c@d_\delta:=c\modcomp e_\delta^{-1},
$
where
\begeq \label{eq:affine-output-feedback-product}
e_\delta=((d_L\circ c)^{\shuffle -1},-(d_L\circ c)^{\shuffle -1}\shuffle (d_R\circ c)),
\endeq
provided that $e_\delta^{-1}\in\allseriesdeltaLC$ and $c@d_\delta\in\allseriesLC$.
In addition, $c@d_{\delta}$ satisfies the fixed-point equation
\begeq \label{eq:affine-output-feedback-fixed-point-equation}
c@d_\delta=c\modcomp (d_L\circ (c@d_\delta),d_R\circ (c@d_\delta)).
\endeq
\endth

\begpr
From Figure~\ref{fig:affine-output-feedback-connection} it is clear that
$u=F_{d_L}[y]v+F_{d_R}[y]=F_{d_L\circ c}[u]v+F_{d_R\circ c}[u]$.
Therefore, using the nonproperness of $d_L\circ c$, and the fact that the shuffle product, shuffle inverse and composition product
on $\allseriesLC$ all preserve local convergence, it follows that
\begin{align*}
u-F_{d_R\circ c}[u]&=F_{d_L\circ c}[u]v \\
F_{(d_L\circ c)^{\shuffle -1}}[u](u-F_{d_R\circ c}[u])&=v \\
F_{((d_L\circ c)^{\shuffle -1},-(d_L\circ c)^{\shuffle -1}\shuffle(d_R\circ c))}[u]&=v.
\end{align*}
It is easily verified that $e:=((d_L\circ c)^{\shuffle -1},-(d_L\circ c)^{\shuffle -1}\shuffle(d_R\circ c))\in \allseriesnp\cap\allseriesdeltaLC$, so if its
inverse is locally convergent, then
$
u=F_{e_\delta^{-1}}[v].
$
In which case, $y=F_c[u]=F_{c\modcomp e_\delta^{-1}}[v]=F_{c@d}[v]$.
Since generating series are always unique, this implies identity \rref{eq:affine-output-feedback-product}.
To determine the fixed-point equation that $c@d_\delta$ must satisfy, observe that the output affine feedback equation
is $y=F_c[F_{d_L}[y]v+F_{d_R}[y]]$. So if $y=F_{c@d_\delta}[v]$ then
\begdi
y=F_c[F_{d_L}[F_{c@d_\delta}[v]]v+F_{d_R}[F_{c@d_\delta}[v]]]=F_{c\modcomp (d_L\circ (c@d_\delta),d_R\circ (c@d_\delta))}[v].
\enddi
Again, since generating series are unique, equation \rref{eq:affine-output-feedback-fixed-point-equation} must hold.
\endpr

\begex \label{ex:cL-equal-one-subgroup}
Consider the feedback system shown in Figure~\ref{fig:affine-output-feedback-connection}
with $d_L=1$. As discussed in the introduction, this corresponds to conventional output feedback and
was the main object of study in
\cite{Gray-etal_SCL14,Gray-Duffaut_Espinosa_CDC13}.
In this case, $F_{c_{\delta}}[u]=u+F_{c_R}[u]$ when $c_R\in\allseriesLC$, and
\rref{eq:affine-output-feedback-product}-\rref{eq:affine-output-feedback-fixed-point-equation} reduce to their expected forms
$
e_{\delta}=(1,-(d_R\circ c))
$
and
$
(c@d_\delta)=c\modcomp (1,d_R\circ(c@d_\delta)),
$
respectively, in Theorem~\ref{th:feedback-product-formula}.
\endex

\begex \label{ex:multiplicative-feedback-formula}
Consider now the case where $d_R=0$ in Figure~\ref{fig:affine-output-feedback-connection}, i.e., pure multiplicative feedback as used, for example,
in phase-locked loops. This problem lies entirely outside the framework of earlier work by the authors.
In this case,
\rref{eq:affine-output-feedback-product}-\rref{eq:affine-output-feedback-fixed-point-equation} reduce to
$
e_{\delta}=((d_L\circ c)^{\shuffle -1},0)
$
and
$
c@d_\delta=c\modcomp (d_L\circ (c@d_\delta),0).
$
respectively.
A specific example exercising these formulas is given in Example~\ref{ex:multiplicative-feedback-formula-computed} after
the Hopf algebra machinery is developed in the next section for computing the composition inverse.
\endex

\section{Hopf Algebra of Coordinate Functions for $\mbf{\allseriesdeltanp}$}
\label{sec:hopf-algebra}

In this section the Hopf algebra of coordinate functions for the group $\allseriesdeltanp$ is described.
This provides an explicit computational framework in which to compute group inverses, and thus, evaluate the output
affine feedback product as described in the previous section.
It is first necessary to restrict the set up to the subset of series having the form $c_{\delta}=(1+c_L^\prime,c_R)$, where $c_L^\prime$ is
proper.
From a control theory point of view, there is no loss of generality since the generating series of any Fliess operator $y=F_{c_\delta}[u]$
can assume this form by simply rescaling $y$. So abusing the notation, in this section $\allseriesdeltanp$ will be used to denote exclusively this subset.
The strategy is to first introduce a connected graded algebra and then a compatible coalgebra in
terms of the composition products.  The connectedness property ensures that this bialgebra is a Hopf algebra.
The section culminates by providing a purely inductive
formula for the antipode of this Hopf algebra, and some illustrative examples are given.

For any word $\eta\in X^\ast$ define the left and right coordinate functions as
\begeq \label{eq:coordinate-maps}
b_\eta:\allseriesdeltanp\rightarrow \re:c_\delta\mapsto (c_L,\eta),\;\;a_\eta:\allseriesdeltanp\rightarrow \re:c_\delta\mapsto (c_R,\eta),
\endeq
respectively.\footnote{This terminology will be justified later.}
Let $V$ denote the $\re$-vector space spanned by these maps.
Define the corresponding free commutative algebra, $H$, with product
$
\mu:h_\eta\otimes \tilde{h}_{\xi}\mapsto h_{\eta}\tilde{h}_{\xi},
$
$h,\tilde{h}\in\{a,b\}$ and
unit $\mbf{1}$ which maps every $c_\delta\in\allseriesnp$ to $1$.
Consequently, definition \rref{eq:coordinate-maps} is extended multiplicatively, that is, for $h,\tilde{h} \in\{a,b\}$ the product $h_{\eta}\tilde{h}_{\xi}$ corresponds to the pointwise product of coordinate functions with respect to $c_\delta=(c_L,c_R) \in \allseriesdeltanp$
\begin{equation}
\label{pointwise}
	h_\eta\tilde{h}_\xi (c_\delta) := h_\eta (c_\delta) \tilde{h}_\xi (c_\delta).
\end{equation}
The {\em degree} of a coordinate function is taken as
$\deg(b_{\eta})=2\abs{\eta}_{x_0}+\abs{\eta}_{x_1}$,
$\deg(a_{\eta})=2\abs{\eta}_{x_0}+\abs{\eta}_{x_1}+1$,
and $\deg(\mbf{1})=0$. In which case,
$V$ is a connected graded vector space, that is, $V=\bigoplus_{n\geq 0} V_n$
with $V_n$ denoting the span of all coordinate functions of degree $n$ and $V_0=\re\mbf{1}$.
Let $V_+=\bigoplus_{n\geq 1} V_n$.
Similarly, $H$ has the connected graduation
$H=\bigoplus_{n\geq 0} H_n$ with $H_0=\re\mbf{1}$.

Three coproducts are now introduced.
The first coproduct is used to define the Hopf algebra on $H$. The remaining two products provide
an inductive method in which to compute it.
Recalling that $(c_\delta\circ d_\delta)_L=(c_L\modcomp d_{\delta})\shuffle d_L$
and $(c_\delta\circ d_\delta)_R=(c_L\modcomp d_{\delta})\shuffle d_R+c_R\modcomp d_{\delta}$,
define coproduct $\Delta$ for any $b_\eta,a_\eta\in V_+$ such that
\begin{subequations} \label{eq:FdB-coproduct-defined}
\begin{align}
\Delta b_{\eta}(c_\delta,d_\delta)&=b_{\eta}(c_\delta\circ d_\delta)=((c_L\modcomp d_{\delta})\shuffle d_L,\eta) \\
\Delta a_{\eta}(c_\delta,d_\delta)&=a_{\eta}(c_\delta\circ d_\delta)=((c_L\modcomp d_{\delta})\shuffle d_R,\eta)+(c_R\modcomp d_{\delta},\eta).
\end{align}
\end{subequations}
The coassociativity of $\Delta$ follows directly from the associativity of the group product on $\allseriesnp$.
The second coproduct is
$\Delta^\cdot_\shuffle(V_+)\subset V_+\otimes V_+$, which is defined for any $h,\tilde{h}\in V_+$ by
\begin{subequations}
\label{eq:shuffle-coproduct-induction}
\begin{align}
\Delta_{\shuffle}^{\tilde{h}}h_{\emptyset}=&h_{\emptyset}\otimes \tilde{h}_{\emptyset} \\
\Delta_{\shuffle}^{\tilde{h}}\circ\theta_k=&(\theta_k\otimes \mbf{1}+\mbf{1}\otimes \theta_k)\circ\Delta_{\shuffle}^{\tilde{h}},
\end{align}
\end{subequations}
where $\theta_k$ denotes the endomorphism on $V_+$ specified by
$\theta_kh_\eta=h_{x_k\eta}$ for $k=0,1$. Clearly this coproduct can be computed recursively.
Finally, consider
for any $\eta\in X^\ast$ the coproduct $\tilde{\Delta}(V_+)\subset V_+\otimes H$, where
\begeq
\tilde{\Delta}b_{\eta}(c_\delta,d_\delta)=(c_L\modcomp d_{\delta},\eta),\;\;
\tilde{\Delta}a_{\eta}(c_\delta,d_\delta)=(c_R\modcomp d_{\delta},\eta). \label{eq:a-and-b-functions}
\endeq
In either case, using
the notation of Sweedler \cite{Sweedler_69},
\begdi
\tilde{\Delta}h_{\eta}(c_\delta,d_\delta)
=\sum h_{\eta(1)}(c_\delta)\;h_{\eta(2)}(d_\delta)
=\sum h_{\eta(1)}\otimes h_{\eta(2)}(c_\delta,d_\delta),
\enddi
where $h_{\eta(1)}\in V_+$ and $h_{\eta(2)}\in H$.
The sums are taken over all the terms that appear in
the respective composition in \rref{eq:a-and-b-functions},
and the specific nature of the factors $h_{\eta(1)}$ and $h_{\eta(2)}$
is not important here.
Like the coproduct $\Delta^\cdot_\shuffle$, this coproduct can also be computed inductively as described next.

\begle \label{le:mixed-copproduct-inductions}
The following identities hold:
\begdes
\item[(1)] $\tilde{\Delta}h_\emptyset=h_{\emptyset}\otimes \mbf{1}$
\item[(2)]
$\tilde{\Delta}\circ \theta_1=
(\theta_1\otimes\mu)\circ(\tilde{\Delta}\otimes \id)\circ\Delta_{\shuffle}^b$
\item[(3)]
$\tilde{\Delta}\circ \theta_0=(\theta_0\otimes \id)\circ \tilde{\Delta}+
(\theta_1\otimes\mu)\circ(\tilde{\Delta}\otimes \id)\circ\Delta_{\shuffle}^a$,
\enddes
where $\id$ denotes the identity map on $H$.
\endle

\begpr
$\;$ \\
\noindent (1) Assume $h=a$ and write $c_R=x_0c_R^0+x_1c_R^1+(c_R,\emptyset)$ with
$c_R^0,c_R^1\in\allseries$. Then from Lemma~\ref{le:mixed-product-properties} (1) and (4) it follows
that
\begin{align*}
c_R\modcomp d_\delta&=(x_0c_R^0)\modcomp d_\delta+(x_1c_R^1)\modcomp d_\delta+(c_R,\emptyset)\modcomp d_\delta \\
&=x_0(c_R^0\modcomp d_\delta)+x_1(d_L\shuffle(c_R^1\modcomp d_\delta))+x_0(d_R\shuffle (c_R^1\modcomp d_\delta))+(c_R,\emptyset).
\end{align*}
In which case, $\tilde{\Delta} a_\emptyset(c_\delta,d_\delta)=(c_R\modcomp d_{\delta},\emptyset)=(c_R,\emptyset)=(a_\emptyset\otimes \mbf{1})(c_\delta,d_\delta)$.
A similar argument holds when $h=b$.

\noindent (2)
If $h=a$ then
\begin{align*}
(\tilde{\Delta}\circ \theta_i) a_\eta(c_\delta,d_\delta)&= \tilde{\Delta}a_{x_i\eta}(c_\delta,d_\delta)=(c_R\modcomp d_\delta,x_i\eta) \\
&=(x_i^{-1}[x_0(c_R^0\modcomp d_\delta)+x_1(d_L\shuffle(c_R^1\modcomp d_\delta))+x_0(d_R\shuffle (c_R^1\modcomp d_\delta))],\eta) \\
&=({\mathds 1}_{i0}[c_R^0\modcomp d_\delta+d_R\shuffle (c_R^1\modcomp d_\delta)]+{\mathds 1}_{i1}[d_L\shuffle(c_R^1\modcomp d_\delta)],\eta),
\end{align*}
where $x_i^{-1}(\cdot)$ is the $\re$-linear left-shift operator specified by
$x_i^{-1}(\eta)=\eta^\prime$ when $\eta=x_i\eta^\prime$ with $\eta^\prime\in X^\ast$ and zero otherwise, and
${\mathds 1}_{xy}$ is the indicator function. So ${\mathds 1}_{xy}=1$ when $x=y$ and zero otherwise.
Letting $c^1_\delta=(c_L,c^1_R)$, it follows that
\begin{align*}
(\tilde{\Delta}\circ \theta_1)a_\eta(c_\delta,d_\delta)&= (d_L\shuffle(c_R^1\modcomp d_\delta),\eta)
= \sum_{\xi,\nu\in X^\ast} (c_R^1\modcomp d_\delta,\xi)(d_L,\nu)(\xi\shuffle \nu,\eta) \\
&= \sum_{\xi,\nu\in X^\ast} \tilde{\Delta}a_\xi(c_\delta^1,d_\delta) b_\nu(d_\delta) (\xi\shuffle \nu,\eta) \\
&= \sum_{\xi,\nu\in X^\ast} \sum a_{\xi(1)}\otimes a_{\xi(2)}(c_\delta^1,d_\delta) b_\nu(d_\delta) (\xi\shuffle \nu,\eta) \\
&= \sum_{\xi,\nu\in X^\ast}\sum \theta_1(a_{\xi(1)})\otimes a_{\xi(2)}(c_\delta,d_\delta) b_\nu(d_\delta) (\xi\shuffle \nu,\eta) \\
&= (\theta_1\circ \id)\circ \sum_{\xi,\nu\in X^\ast} \tilde{\Delta}a_\xi\otimes b_\nu(c_\delta,d_\delta, d_\delta) (\xi\shuffle \nu,\eta) \\
&=(\theta_1\otimes\mu)\circ(\tilde{\Delta}\otimes \id)\circ \Delta_\shuffle^b a_\eta(c_\delta,d_\delta).
\end{align*}
The proof when $h=b$ is perfectly analogous.

\noindent (3)
If $h=a$ then $(\tilde{\Delta}\circ \theta_0) a_\eta(c_\delta,d_\delta)
=(c_R^0\modcomp d_\delta,\eta)+(d_R\shuffle (c_R^1\modcomp d_\delta),\eta)$.
From this point, the method of proof is exactly the same as that in part (2) modulo the fact that $\Delta_\shuffle^a$ is used here
due to the presence of $d_R$ instead of $d_L$ in the shuffle product.
\endpr

\begex
Applying the identities in Lemma~\ref{le:mixed-copproduct-inductions} gives the first few $\tilde{\Delta}$ coproduct terms ordered by degree
$n_b$ $(h=b)$, $n_a$ $(h=a)$:
\begin{align*}
n_b,n_a=0,1 &:\tilde{\Delta}h_{\emptyset}=h_\emptyset \otimes \mbf{1} \\
n_b,n_a=1,2 &:\tilde{\Delta}h_{x_1}= h_{x_1}\otimes \mbf{1} \\
n_b,n_a=2,3 &:\tilde{\Delta}h_{x_0}= h_{x_0}\otimes \mbf{1}+h_{x_1}\otimes a_{\emptyset} \\
n_b,n_a=2,3 &:\tilde{\Delta}h_{x_1^2}= h_{x_1^2}\otimes \mbf{1}+h_{x_1}\otimes b_{x_1} \\
n_b,n_a=3,4 &:\tilde{\Delta}h_{x_0x_1}=h_{x_0x_1}\otimes \mbf{1}+h_{x_1}\otimes a_{x_1}+h_{x_1^2}\otimes a_{\emptyset} \\
n_b,n_a=3,4 &:\tilde{\Delta}h_{x_1x_0}=h_{x_1x_0}\otimes \mbf{1}+h_{x_1}\otimes b_{x_0}+h_{x_1^2}\otimes a_{\emptyset} \\
n_b,n_a=3,4 &:\tilde{\Delta}h_{x_1^3}=
h_{x_1^3}\otimes \mbf{1} + 3 h_{x_1^2}\otimes b_{x_1}+h_{x_1}\otimes b_{x_1^2} \\
n_b,n_a=4,5 &:\tilde{\Delta}h_{x_0^2}=h_{x_0^2}\otimes \mbf{1}+h_{x_1} \otimes a_{x_0}+h_{x_0x_1}\otimes a_{\emptyset}+
h_{x_1x_0}\otimes a_{\emptyset}+h_{x_1^2}\otimes (a_{\emptyset})^2.
\end{align*}

\endex

The next lemma provides a grading for the coproduct $\tilde{\Delta}$, which is clearly evident in the example above.
Its particular form, namely having $V$ on the left,
implies that $H$ is a right-sided Hopf algebra \cite{Menous-Patras_16,Loday-Ronco_10}, a property that is essential
in Section \ref{sec:FPS-Lie-group}.

\begle \label{le:modified-coproduct-grading}
For any $h_\eta\in V_n$
\begeq \label{eq:grading-mixed-coproduct}
\tilde{\Delta}h_\eta\in \bigoplus_{j+k=n} V_j\otimes H_k=:(V\otimes H)_n.
\endeq
\endle
\begpr
The following facts are essential:
\begdes
\item[(1)] $\deg(\theta_1h)=\deg(h)+1$
\item[(2)] $\deg(\theta_0h)=\deg(h)+2$
\item[(3)] $\Delta_\shuffle^{\tilde{h}} h \in (V\otimes V)_{\deg(h)+{\mathds 1}_{\tilde{h}a}}$.
\enddes
The proof is via
induction on the length of $\eta$. When $\abs{\eta}=0$ then clearly
$\tilde{\Delta}h_\emptyset=h_\emptyset\otimes \mbf{1}\in V_n\otimes H_0$,
where $n=\deg(h_\emptyset)\in\{0,1\}$ (noting that $b_\emptyword\sim\mbf{1}$).
Assume now that \rref{eq:grading-mixed-coproduct} holds for words up to length
$\abs{\eta}\geq 0$. Let $n=\deg(h_\eta)$.
There are two ways to increase the
length of $\eta$.
First consider $h_{x_1\eta}$.
From item (1) above $\deg(h_{x_1\eta})=n+1$. Now apply item (3), the induction hypothesis, and Lemma~\ref{le:mixed-copproduct-inductions} in that order:
\begin{align*}
\Delta_{\shuffle}^b h_\eta & \in (V\otimes V)_n \\
(\tilde{\Delta}\otimes \id)\circ \Delta_{\shuffle}^b h_\eta & \in
(V\otimes H\otimes V)_n \\
(\theta_1\otimes\mu)\circ(\tilde{\Delta}\otimes \id)\circ \Delta_{\shuffle}^b h_\eta & \in \bigoplus_{j+k=1}^n\in V_{j+1}\otimes H_k \\
\tilde{\Delta}h_{x_1\eta}&\in(V\otimes H)_{n+1},
\end{align*}
which proves the assertion.
Consider next $h_{x_0\eta}$. From item~(2) above $\deg(h_{x_0\eta})=n+2$. In this case, repeat the first two steps of the previous case and apply item (1) to get
\begdi
(\theta_1\otimes\mu)\circ(\tilde{\Delta}\otimes \id)\circ \Delta_{\shuffle}^a h_\eta  \in \bigoplus_{j+k=1}^{n+1}\in V_{j+1}\otimes H_k\subset (V\otimes H)_{n+2}.
\enddi
In addition, from the induction hypothesis and item~(2) it follows that
\begdi
(\theta_0\otimes\id)\circ \tilde{\Delta}h_\eta\in\bigoplus_{j+k=1}^n V_{j+2}\otimes H_k\subset (V\otimes H)_{n+2}.
\enddi
Thus, applying Lemma~\ref{le:mixed-copproduct-inductions}, $\tilde{\Delta}h_{x_0\eta}\in(V\otimes H)_{n+2}$,
which again proves the assertion
and completes the proof.
\endpr

The following lemma is of central importance as it shows how to compute the Hopf algebra coproduct $\Delta$ in terms
of $\tilde{\Delta}$ and $\Delta_\shuffle$.

\begle \label{le:general-group-copproduct-inductions}
The following identities holds:
\begdes
\item[(1)] $\Delta b_\eta=(\id\otimes\mu)\circ(\tilde{\Delta}\otimes\id)\circ \Delta^b_{\shuffle} b_\eta$
\item[(2)] $\Delta a_\eta=(\id\otimes\mu)\circ(\tilde{\Delta}\otimes\id)\circ \Delta^a_{\shuffle} b_\eta+\tilde{\Delta}a_\eta$.
\enddes
\endle

\begpr

\noindent (1) Observe
\begin{align*}
\Delta b_{\eta}(c_\delta,d_\delta)&= \sum_{\xi,\nu\in X^\ast} (c_L\modcomp d_\delta,\xi)(d_L,\nu)(\xi\shuffle\nu,\eta)
=\sum_{\xi,\nu\in X^\ast} \tilde{\Delta}b_{\xi}(c_\delta,d_\delta)b_{\nu}(d_\delta)(\xi\shuffle\nu,\eta) \\
&=(\tilde{\Delta}\otimes\id)\circ \Delta_\shuffle^b b_{\eta}(c_\delta,d_\delta,d_\delta)
=(\id\otimes\mu)\circ(\tilde{\Delta}\otimes\id)\circ \Delta^b_{\shuffle} b_\eta(c_\delta,d_\delta).
\end{align*}
(2) In a similar fashion
\begin{align*}
\Delta a_{\eta}(c_\delta,d_\delta)
&= (\tilde{\Delta}\otimes \id)\circ \Delta_\shuffle^a b_\eta(c_\delta,d_\delta,d_\delta) + \tilde{\Delta} a_\eta (c_\delta,d_\delta) \\
&= [(\id\otimes \mu) \circ (\tilde{\Delta}\otimes \id)\circ \Delta_\shuffle^a b_\eta + \tilde{\Delta} a_\eta] (c_\delta,d_\delta).
\end{align*}
\endpr

\begex
Applying the identities in Lemma~\ref{le:general-group-copproduct-inductions} gives the first few reduced coproduct terms, namely,
$\Delta^\prime h_\eta:=\Delta h_\eta- h_\eta\otimes \mbf{1} - \mbf{1}\otimes h_\eta$:
\begin{align*}
n=1 &: \Delta^\prime b_{x_1}= 0 \\
n=2 &: \Delta^\prime b_{x_0}= b_{x_1} \otimes a_{\emptyset} \\
n=2 &: \Delta^\prime b_{x_1^2}= 3 b_{x_1} \otimes b_{x_1} \\
n=3 &: \Delta^\prime b_{x_0x_1}= b_{x_0} \otimes b_{x_1}  + b_{x_1} \otimes b_{x_0} +
b_{x_1} \otimes  a_{x_1} +
b_{x_1} \otimes b_{x_1} a_{\emptyset} + b_{x_1^2} \otimes a_{\emptyset} \\
n=3 &: \Delta^\prime b_{x_1x_0}=b_{x_0} \otimes b_{x_1} + 2 b_{x_1} \otimes b_{x_0} +
b_{x_1} \otimes b_{x_1} a_{\emptyset} + b_{x_1^2} \otimes a_{\emptyset} \\
n=3 &: \Delta^\prime b_{x_1^3}=6b_{x_1^2}\otimes b_{x_1}+4b_{x_1}\otimes b_{x_1^2}+  3b_{x_1}\otimes (b_{x_1})^2\\
n=4 &: \Delta^\prime b_{x_0^2}=2 b_{x_0} \otimes b_{x_0} + b_{x_1} \otimes a_{x_0} +
2 b_{x_1} \otimes b_{x_0} a_{\emptyset} +
b_{x_0x_1} \otimes a_{\emptyset} + b_{x_1x_0} \otimes a_{\emptyset} + \\
&\hspace*{0.18in}b_{x_1^2} \otimes (a_{\emptyset})^2
\end{align*}
\begin{align*}
n=1 &: \Delta^\prime a_{\emptyset}= 0 \\
n=2 &: \Delta^\prime a_{x_1}= b_{x_1} \otimes  a_{\emptyset} \\
n=3 &: \Delta^\prime a_{x_0}=  b_{x_0} \otimes a_{\emptyset} + a_{x_1} \otimes a_{\emptyset} + b_{x_1} \otimes (a_{\emptyset})^2 \\
n=3 &: \Delta^\prime a_{x_1^2}= a_{x_1} \otimes b_{x_1} + 2 b_{x_1} \otimes a_{x_1} +
b_{x_1} \otimes b_{x_1} a_{\emptyset} +
 b_{x_1^2} \otimes a_{\emptyset} \\
n=4 &: \Delta^\prime a_{x_0x_1}= b_{x_1} \otimes a_{x_0} + b_{x_0} \otimes a_{x_1} +
 a_{x_1} \otimes a_{x_1} + b_{x_0x_1} \otimes a_{\emptyset} +
a_{x_1^2} \otimes a_{\emptyset} +  \\
&\hspace*{0.18in}2 b_{x_1} \otimes a_{x_1} a_{\emptyset} + b_{x_1^2} \otimes (a_{\emptyset})^2 \\
n=4 &: \Delta^\prime a_{x_1x_0}= a_{x_1} \otimes b_{x_0} + b_{x_1} \otimes a_{x_0} + b_{x_0}\otimes a_{x_1} +
b_{x_1}\otimes b_{x_0} a_{\emptyset} +
b_{x_1x_0} \otimes a_{\emptyset} +  \\
&\hspace*{0.18in} a_{x_1^2} \otimes a_{\emptyset} + b_{x_1} \otimes a_{x_1} a_{\emptyset} + b_{x_1^2} \otimes (a_{\emptyset})^2 \\
n=4 &: \Delta^\prime a_{x_1^3}=b_{x_1^3}\otimes a_\emptyset + 3b_{x_1^2}\otimes b_{x_1}a_\emptyset + b_{x_1}\otimes b_{x_1^2}a_\emptyset +
3b_{x_1^2}\otimes a_{x_1} +3 b_{x_1}\otimes b_{x_1}a_{x_1}+ \\
&\hspace*{0.18in} 3b_{x_1}\otimes a_{x_1^2}+
3a_{x_1^2}\otimes b_{x_1} + a_{x_1}\otimes b_{x_1^2} \\
n=5 &: \Delta^\prime a_{x_0^2}=  2 b_{x_0} \otimes a_{x_0} + a_{x_1} \otimes a_{x_0} +
b_{x_0^2} \otimes a_{\emptyset} +  a_{x_0x_1} \otimes a_{\emptyset} +
a_{x_1x_0} \otimes a_{\emptyset} +   \\
&\hspace*{0.18in} 3 b_{x_1} \otimes a_{x_0} a_{\emptyset} + b_{x_0x_1} \otimes (a_{\emptyset})^2 +
b_{x_1x_0} \otimes (a_{\emptyset})^2 +
a_{x_1^2} \otimes (a_{\emptyset})^2 +
b_{x_1^2} \otimes (a_{\emptyset})^3.
\end{align*}
\endex

The main result of this section is given next, namely that the commutative unital algebra $(H,\mu)$ together with the coproduct $\Delta$ is a connected graded bialgebra, and thus has an antipode $S$. Before stating the actual theorem, some additional motivation is given concerning the nature of the antipode. Observe that for any $c_\delta \in \allseriesdeltanp$ one can
naturally associate a Hopf algebra character $\Phi_c \in L(H,\re)$ by letting
$\Phi_c:b_\eta\mapsto (c_L,\eta)$, $\Phi_c:a_\eta\mapsto (c_R,\eta)$,
and $\Phi_c(\mbf{1})=1$. Note that \rref{pointwise} implies $\Phi_c (h_\eta \tilde{h}_\xi)=\Phi_c(h_\eta)\Phi_c(\tilde{h}_\xi)$ for any $h,\tilde{h}\in \{a,b\}$ and $\eta,\xi \in X^\ast$. Since the coproduct $\Delta$ of $H$ has been constructed to satisfy \rref{eq:FdB-coproduct-defined}, the following simple calculation shows that for all $c_\delta, d_\delta \in \allseriesdeltanp$, the convolution product of $\Phi_{c}, \Phi_{d} \in L(H,\re)$ corresponds to composition in $\allseriesdeltanp$:
\begdi
\label{eq:general-coordinate-maps}
(\Phi_{c} \star \Phi_{d})(h_\eta) := m_\re(\Phi_c \otimes \Phi_d)\Delta(h_\eta)=h_\eta(c_\delta \circ d_\delta)=(c_\delta \circ d_\delta,\eta).
\enddi%
(See \cite[Lemma 2]{Gray-etal_SCL14} for a similar calculation.) In light of \rref{eq:convinv-via-antipode}, it then follows that
\begeq \label{eq:c-delta-inverse-via-S}
h_\eta(c_\delta^{-1})=(S h_\eta)(c_\delta),\;\;\forall \eta\in X^\ast,\;\; h\in \{b,a\}.
\endeq
Thus, the antipode provides an explicit way to compute the group inverse.

\begth \label{th:hopf-algebra-for-coordinate-functions}
$(H,\mu,\Delta)$ is a connected graded commutative unital Hopf algebra.
\endth

\begpr
From the development above, it is clear that $(H,\mu,\Delta)$
is a connected bialgebra with unit $\mbf{1}$ and a counit $\varepsilon$ defined by
$\varepsilon(a_\eta)=0$ for all $\eta\in X^\ast$, $\varepsilon(b_\eta)=0$ for all nonempty
$\eta\in X^\ast$, and
$\varepsilon(\mbf{1})=1$.
Here it is shown that this bialgebra is also graded
and thus is automatically a Hopf algebra, i.e.,
has a well defined antipode, $S$ \cite{Figueroa-Gracia-Bondia_05}.
It only remains to be shown for any $n\geq 0$
that
$
\Delta H_n\subseteq (H\otimes H)_n.
$
It is well known if $h\in V_n$ then $\Delta_\shuffle^{\tilde{h}}h\in (V\otimes V)_n$.
Therefore, it follows directly from Lemmas~\ref{le:modified-coproduct-grading} and
\ref{le:general-group-copproduct-inductions} that
$
\Delta h\in (V\otimes H)_{n}.
$
In which case, via the identity $\Delta(a^i_\eta a^j_\xi)=\Delta a_\eta^i \Delta a_\xi^j$,
it must hold that
$\Delta H_n\subseteq (H\otimes H)_n$, $n\geq 0$.
\endpr

The next theorem says that the antipode of any connected graded Hopf algebra can be computed in a
recursive manner once the coproduct is computed.

\begth {\rm \cite{Figueroa-Gracia-Bondia_05}} \label{th:antipode-induction}
The antipode, $S$, of any connected graded Hopf algebra $(H,\mu,\Delta)$
can be computed for any $a\in H_{k}$, $k\geq 1$ by
\begdi
S a=-a-\sum (S a^\prime_{(1)})a^\prime_{(2)}=-a-\sum  a^\prime_{(1)}S a^\prime_{(2)},
\label{eq:general-antipode-induction}
\enddi
where the reduced coproduct is
$\Delta^\prime a=\Delta a-a\otimes \mbf{1}-\mbf{1}\otimes a=\sum a^\prime_{(1)}a^\prime_{(2)}$.
\endth

As noted earlier,
the coproducts $\Delta_\shuffle$ and $\tilde{\Delta}$ can be computed recursively, and $\Delta$ is computed directly in terms of
$\Delta_\shuffle$ and $\tilde{\Delta}$. So in fact the
antipode of $H$ can be computed in a {\em fully} recursive manner as described next.

\begth \label{th:fully-recursive-S}
The antipode, $S$, of any $h_\eta\in V_+$ can be computed by the following algorithm:
\begdes
\item[i.] Recursively compute $\Delta_{\shuffle}^{\tilde{h}}$ via \rref{eq:shuffle-coproduct-induction}.
\item[ii.] Recursively compute $\tilde{\Delta}$ via Lemma~\ref{le:mixed-copproduct-inductions}.
\item[iii.] Compute $\Delta$ via Lemma~\ref{le:general-group-copproduct-inductions}.
\item[iv.] Recursively compute $S$ via Theorem~\ref{th:antipode-induction}.
\enddes
\endth

\begex
The first few antipode terms computed via Theorem~\ref{th:fully-recursive-S} are:
\begin{align*}
n=1&: Sb_{x_1}= -b_{x_1}\\
n=2&: Sb_{x_0}= -b_{x_0}+b_{x_1}a_\emptyset \\
n=2&: Sb_{x_1^2}= -b_{x_1^2}+3(b_{x_1})^2 \\
n=3&: Sb_{x_0x_1}= -b_{x_0x_1}+b_{x_1^2}a_{\emptyset}-3 (b_{x_1})^2 a_\emptyset+2 b_{x_0}b_{x_1}+b_{x_1}a_{x_1} \\
n=3&: Sb_{x_1x_0}= -b_{x_1x_0}+b_{x_1^2} a_{\emptyset}+3 b_{x_0}b_{x_1}-3 (b_{x_1})^2 a_{\emptyset} \\
n=3&: Sb_{x_1^3}=-b_{x_1^3}+10b_{x_1}b_{x_1^2}-15(b_{x_1})^3 \\
n=4&: Sb_{x_0^2}= -b_{x_0^2}+b_{x_0x_1} a_\emptyset+b_{x_1x_0}a_\emptyset- b_{x_1^2} (a_\emptyset)^2
+3 (b_{x_1})^2 (a_\emptyset)^2+2 (b_{x_0})^2-  \\
&\hspace*{0.18in} 5 b_{x_0}b_{x_1}a_{\emptyset}-b_{x_1} a_{\emptyset}a_{x_1}+b_{x_1}a_{x_0}
\end{align*}
\begin{align*}
n=1&: Sa_{\emptyset}=-a_\emptyset \\
n=2&: Sa_{x_1}= -a_{x_1}+b_{x_1} a_\emptyset \\
n=3&: Sa_{x_0}= -a_{x_0}+b_{x_0} a_\emptyset- b_{x_1} (a_\emptyset)^2+a_\emptyset a_{x_1} \\
n=3&: Sa_{x_1^2}= -a_{x_1^2}+b_{x_1^2} a_\emptyset - 3 (b_{x_1})^2 a_\emptyset+3 b_{x_1}a_{x_1} \\
n=4&: Sa_{x_0x_1}= -a_{x_0x_1}+b_{x_0x_1} a_{\emptyset}-2 b_{x_0}b_{x_1}a_\emptyset
-4 b_{x_1}a_{\emptyset}a_{x_1}+3 (b_{x_1})^2 (a_{\emptyset})^2- \\
&\hspace*{0.18in} b_{x_1^2} (a_{\emptyset})^2+(a_{x_1})^2+a_\emptyset a_{x_1^2}+b_{x_0}a_{x_1}+b_{x_1} a_{x_0} \\
n=4&: Sa_{x_1x_0}=-a_{x_1x_0}+b_{x_1x_0}a_\emptyset+a_\emptyset a_{x_1^2}-3 b_{x_0}b_{x_1} a_{\emptyset}
-3 b_{x_1}a_\emptyset a_{x_1}-b_{x_1^2} (a_\emptyset)^2+ \\
&\hspace*{0.18in} 3 (b_{x_1})^2(a_\emptyset)^2+2b_{x_0}a_{x_1}+ b_{x_1}a_{x_0} \\
n=4&: Sa_{x_1^3}=-a_{x_1^3}+b_{x_1^3}a_{\emptyset}-10b_{x_1}b_{x_1^2}a_\emptyset +
15(b_{x_1})^3a_\emptyset+4b_{x_1^2}a_{x_1}-15(b_{x_1})^2a_{x_1}+ \\
&\hspace*{0.18in} 6b_{x_1}a_{x_1^2} \\
n=5&: Sa_{x_0^2}=-a_{x_0^2} + b_{x_0^2} a_\emptyset-2 (b_{x_0})^2 a_{\emptyset}+ a_\emptyset a_{x_0x_1}+a_{\emptyset}a_{x_1x_0}
-3 b_{x_1}a_{x_0} a_\emptyset -  \\
&\hspace*{0.18in} 3 b_{x_0}a_{x_1} a_\emptyset -a_\emptyset (a_{x_1})^2-b_{x_0x_1}(a_\emptyset)^2 - b_{x_1x_0}(a_{\emptyset})^2
+5 b_{x_0}b_{x_1} (a_\emptyset)^2-(a_\emptyset)^2a_{x_1^2}+\\
&\hspace*{0.18in} 4 b_{x_1} (a_\emptyset)^2 a_{x_1}+b_{x_1^2} (a_\emptyset)^3-3 (b_{x_1})^2 (a_{\emptyset})^3+2 b_{x_0} a_{x_0} + a_{x_0}a_{x_1}.
\end{align*}
\endex

\begex \label{ex:multiplicative-feedback-formula-computed}
Reconsider the multiplicative feedback system described in Example~\ref{ex:multiplicative-feedback-formula}. Suppose
$c=\sum_{k\geq 0} k!\,x_1^k=\sum_{k\geq 0} x_1^{\shuffle k}=(1-x_1)^{\shuffle -1}$ and $F_{d_\delta}=I$ (unity feedback).
In this case, \rref{eq:affine-output-feedback-fixed-point-equation} reduces to $(c@d_{\delta})=c\modcomp (d_L\circ c@d_\delta,0)$ so that
\begdi
c@d_\delta=c\modcomp{(c^{\shuffle -1},0)^{-1}}=(c^{\shuffle -1},0)^{-1}=(1-x_1,0)^{-1}.
\enddi
The composition inverse can be computed directly from \rref{eq:c-delta-inverse-via-S} and the antipode formulas above.
Namely, taking into account only the antipode terms containing powers of the coordinate function $b_{x_1}$ (all others will yield zero),
it follows that
\begin{align*}
c@d_\delta&=(1-x_1,0)^{-1}
=\left(\left[\mathbf{1}+\sum_{\eta\in X^+} \eta Sb_\eta\right](1-x_1),\sum_{\eta\in X^\ast} \eta Sa_\eta(0)\right) \\
&=([\mathbf{1}+x_1Sb_{x_1}+x_1^2 Sb_{x_1^2}+x_1^3 Sb_{x_1^3}+\cdots](1-x_1),0) \\
&=(1+x_1(-b_{x_1}(1-x_1))+x_1^2([-b_{x_1^2}+3(b_{x_1})^2](1-x_1)+ \\
&\hspace*{0.2in}x_1^3([b_{x_1^3}+10b_{x_1}b_{x_1^2}-15(b_{x_1})^3](1-x_1)),0) \\
&=(1+x_1+3x_1^2+15x_1^3+\cdots,0).
\end{align*}
\endex

The deferred proof from Section~\ref{sec:transformation-group-allseriesdelta} is presented next.

\begprx{of Lemma~\ref{le:mixed-product-properties} (3)}
The assertion is that $c\modcomp d_\delta=k$ implies $c=k$.
The proof is by induction on the grading of $H$.
If $c\modcomp d_\delta=k$ then clearly
$k=a_{\emptyset}(c_\delta\modcomp d_\delta)=\tilde{\Delta}a_{\emptyset}(c_\delta,d_\delta)=a_{\emptyset}(c_{\delta})$ assuming
without loss of generality that $c_\delta=(1,c)$.  Therefore, $(c,\emptyset)=k$.
Similarly, it follows that $0=a_{x_1}(c_\delta\modcomp d_\delta)=\tilde{\Delta}a_{x_1}(c_\delta,d_\delta)=a_{x_1}(c_\delta)$.
Thus, $(c,x_1)=0$.
Now suppose $a_\eta(c_\delta)=0$ for all $a_{\eta}\in H_n$ up to some fixed $n\geq 2$.
Then for any $x_j\in X$
\begdi
0=\tilde{\Delta}a_{x_j\eta}(c_\delta,d_\delta)=a_{x_j\eta}(c_\delta)+
\sum_{a_{x_j\eta(2)}\neq 1} a_{x_j\eta(1)}(c_\delta)\;a_{x_j\eta(2)}(d_\delta),
\enddi
where in general $a_{x_j\eta(1)}\neq a_{\emptyset}$. Therefore,
$a_{x_j\eta}(c_\delta)=0$, or equivalently, $(c,x_j\eta)=0$. In which, case $c=k$.
\endprx

The section is concluded by some dimensional analysis of the grading of $V$ and $H$.
This information could be useful in determining the complexity of the
antipode recursion as a function of the degree of the coordinate functions,
but that topic is beyond the scope of this paper.
Let $V_{h,k}$ denote the subspace of $V_k$ spanned by the coordinate functions $h_\eta$
of degree $k$ where $h\in\{a,b\}$.
Define $p_{h,k}=\dim(V_{h,k})$, $p_k=\dim(V_k)$ and the corresponding generating functions
$F_{V_h}=\sum_{k\geq 1} p_{h,k}X^k$, $F_V=\sum_{k\geq 1} p_k X^k$. Analogous definitions apply
when $V$ is replaced by $H$.
\begth
The following identities hold:
\begin{align*}
F_{V_a}&= \frac{X}{1-X-X^2}
= X+X^2+2 X^3+3 X^4+5 X^5+8 X^6+13 X^7+21 X^8+\cdots \\
F_{V_b}&= \frac{X+X^2}{1-X-X^2}
= X+2 X^2+3 X^3+5 X^4+8 X^5+13 X^6+21 X^7+34 X^8+\cdots \\
F_V&= F_{V_a}+F_{V_b} =\frac{2X+X^2}{1-X-X^2}
=2 X+3 X^2+5 X^3+8 X^4+13 X^5+21 X^6+\cdots \\
F_{H_a}& =\prod_{k=1}^\infty \frac{1}{(1-X^k)^{p_{a,k}}}
=1+X+2 X^2+4 X^3+8 X^4+15 X^5+30 X^6+56 X^7+\cdots \\
F_{H_b}&= \prod_{k=1}^\infty \frac{1}{(1-X^k)^{p_{b,k}}}
=1+X+3 X^2+6 X^3+14 X^4+28 X^5+61 X^6+122 X^7+\cdots \\
F_{H}&= \prod_{k=1}^\infty \frac{1}{(1-X^k)^{p_{k}}}=F_{H_a}F_{H_b}
= 1+2 X+6 X^2+15 X^3+38 X^4+89 X^5+\cdots.
\end{align*}
\endth
\begpr
The identity for $F_{V_a}$ is proved in \cite[Proposition 8]{Foissy_CiA15}, the proof for $F_{V_b}$ is perfectly analogous.
The identity for $F_V$ follows directly from the fact that $V=V_a \oplus V_b$. It is worth noting that
the coefficients of all three series come from the Fibonacci sequence.
The identity for $F_{H_a}$ was also proved in \cite{Foissy_CiA15}, and again the proof for $F_{H_b}$ is very similar.
The factorization of $F_H$ is a consequence of the fact that $p_k=p_{a,k}+p_{b,k}$. In this case, the coefficients of
$F_{H_a}$ and $F_{H_b}$ are integer sequences A166861 and A200544, respectively, in \cite{OEIS}, while the sequence
for $F_H$ appears to be new.
\endpr

\section{System Inversion and Feedback Invariants}
\label{sec:group-invariants}

In this section two system theoretic problems connected with the output affine feedback
transformation group are addressed. First it is shown how to pose the dynamic inversion in terms of this transformation
group and its inverse. This naturally leads to the topic of feedback input-output linearization
via nonlinear feedback.
This then takes the presentation to the second problem of feedback invariants under the output affine feedback
transformation group.

\subsection{System Inversion}
A problem often encountered in control applications is the following system inverse
problem. An input $u$ is mapped to an output $y$ by an operator $F$ defined in terms of a dynamical system
\begin{subequations} \label{eq:siso-space-realization}
\begin{align}
\dot{z}&= g_0(z)+g_1(z)u,\;\; z(t_0)=z_0 \\
y&=h_L(z)u+h_R(z),
\end{align}
\end{subequations}
where all the functions of
the state $z$ are assumed to be analytic on some neighborhood ${\cal W}$ of $z_0$.
If $h_{L}(z_0)\neq 0$ then $F$ is locally
invertible. The left inverse, $F^{-1}:y\mapsto u$, can be determined dynamically by first solving
the output equation for $u$, $u=h_L^{-1}(z)[y-h_R(z)]$, and
then substituting the result into the state equation so that the system
\begin{align*}
	\dot{z}&= \left[g_0(z)-g_1(z)h^{-1}_{L}(z)h_{R}(z)\right]+g_1(z) h_{L}^{-1}(z)y,\;\; z(t_0)=z_0 \\
u&=h_L^{-1}(z)y-h_L^{-1}(z)h_R(z)
\end{align*}
has the property that if $y$ is a given function in the range of $F$ then $y=F[u]$ maps to $u$
over an interval $[t_0,t_0+T]$ with $T>0$.
In short, input-output inversion is done by mapping the given system $(g_0,g_1,h_L,h_R,z_0)$ to its inverse system
\begin{equation} \label{eq:realization2inverse-realization}
	(\bar{g}_0,\bar{g}_1,\bar{h}_L,\bar{h}_R,z_0)
	:=\left(g_0-g_1h_{L}^{-1}h_{R},g_1 h^{-1}_{L},h_L^{-1},-h_L^{-1}h_R,z_0\right).
\end{equation}
It is not difficult to see that the mapping $F:u\mapsto y$ can be described independent of the state
$z$ and its dynamics by an element $c_\delta=(c_L,c_R)$ of the output affine feedback transformation group so that
$y=F[u]=F_{c_L}[u]u+F_{c_R}[u]=:F_{c_\delta}[u]$, where
$c_L,c_R \in\allseriesLC$ are each computed from \rref{eq:c-equals-Lgh} using the respective outputs $h_L$ and $h_R$.
The condition that $h_L(z_0)\neq 0$ ensures that $c_L$ is nonproper. In which case, the inverse input-output map
is given by $F_{c_\delta^{-1}}$, where the generating series $c_\delta^{-1}$ can be computed using the antipode of the underlying Hopf algebra
described in Section~\ref{sec:hopf-algebra}. In fact, the following example demonstrates that this relation goes in both directions when
$c_\delta=(c_L,c_R)$ has a realization of the form \rref{eq:siso-space-realization}.

\begex
Reconsider Example~\ref{ex:multiplicative-feedback-formula-computed}
where it was necessary to compute $(1-x_1,0)^{-1}$.
It is straightforward to show that $y_L=F_{c_L}[u]=F_{1-x_1}[u]$, $y_R=0$ is realized by
$(g_0,g_1,h_L,h_R,z_0)=(0,-1,z,0,1)$ so that
by \rref{eq:realization2inverse-realization} it follows that
$(\bar{g}_0,\bar{g}_1,\bar{h}_L,\bar{h}_R,z_0)=(0,-1/z,1/z,0,1)$. In which
case, $c_R^{\circ -1}=0$ and
\begin{align*}
c_L^{\circ -1}&=1+x_1L_{\bar{g}_1}\bar{h}_L(z_0)+x_1^2L^2_{\bar{g}_1}\bar{h}_L(z_0)+x_1^3L^3_{\bar{g}_1}\bar{h}_L(z_0)+\cdots \\
&=1+x_1+3x_1^2+15x_1^3+\cdots,
\end{align*}
as computed earlier.
\endex

Given the connection between dynamic inversion and
feedback linearization \cite{Isidori_95,Nijmeijer-vanderSchaft_90}, it is natural to consider the latter
in terms of the output affine feedback transformation group.
Consider a Fliess operator
$y=F_c[u]$ which has an $n$ dimensional analytic state space realization of the form
\begeq \label{eq:control-affine-realization}
\dot{z}=g_0(z)+g_1(z)u,\;\;z(0)=z_0,\;\; y=h(z),
\endeq
with relative degree $r$. In which case, it follows that
$
y^{(r)}=L_{g_0}^rh(z)+L_{g_1} L_{g_0}^{r-1}h(z)u
$
with $L_{g_1}L_{g_0}^{r-1}h(z)$ being nonzero on some neighborhood of $z_0$ \cite{Isidori_95,Nijmeijer-vanderSchaft_90}.
Since the state can be viewed in terms of a Fliess operator for some generating series $c_z\in\allseriesnLC$, it follows
from Fliess's fundamental formula that there are generating series $e_L,e_R\in\allseriesLC$ such that $F_{e_R}[u]=L_{g_0}^rh(F_{c_z}[u])$
and $F_{e_L}[u]=L_{g_1}L_{g_0}^{r-1}h(F_{c_z}[u])$, respectively, so that
$
y^{(r)}=F_{e_L}[u]u+F_{e_R}[u].
$
Since $(e_L,\emptyset)=L_{g_1}L_{g_0}^{r-1}h(F_{c_z}[u])(0)\neq 0$, it follows that $e_\delta=(e_L,e_R)\in\allseriesdeltanp\cap\allseriesdeltaLC$, and therefore, if $v:=y^{(r)}$
then the feedback linearization law can be written in the form $u=F_{e_\delta^{-1}}[v]$. When this feedback is applied to the plant, $F_c$,
the feedback linearized system is described by
$$
	y=F_c[u]=F_c[F_{e_\delta^{-1}}[v]]=F_{c\modcomp e_\delta^{-1}}[v]=v^{(-r)}.
$$
So it is evident that the
output affine feedback transformation group is at play in this problem. The proposition, however, is that this analysis holds even when $F_c$ does {\em not} have a
state space realization.

The starting point for a realization free framework for feedback input-output linearization is the
following definition which describes relative degree from a
generating series
point of view. It uses the notion of
a {\em linear word}, that is, any word in the language
$
{\mathscr L}=\{\eta\in X^{\ast}:\eta=x_0^{n_1}x_1x_0^{n_0},\;n_1,n_0\geq 0\}.
$
Furthermore, note that every $c\in\allseries$ can be decomposed into its natural and forced components, that is,
$c=c_N+c_F$, where $c_N:=\sum_{k\geq 0} (c,x_0^k)x_0^k$
and $c_F:=c-c_N$.
Finally, for any letter $x_i\in X$, the left-shift operator $x_i^{-1}(\cdot)$ introduced
earlier is defined
inductively for higher order shifts via $(x_i\xi)^{-1}(\cdot)=\xi^{-1}x_i^{-1}(\cdot)$,
where $\xi\in X^\ast$.

\begde\cite{Gray-etal_Auto14} \label{de:relative-degree-c}
Given $c\in\allseries$, let $r\geq 1$ be the largest integer such that $\supp(c_F)\subseteq x_0^{r-1}X^\ast$.
Then $c$ has \bfem{relative degree} $r$ if the linear word $x_0^{r-1}x_1\in \supp(c)$,
otherwise it is not well defined.
\endde

Observe that $c$ having relative degree $r$ is equivalent to saying that
\begeq \label{eq:cF-relative-degree-decomposition}
c=c_N+c_F=c_N+Kx_0^{r-1}x_1+x_0^{r-1}e
\endeq
for some $K\neq 0$ and some proper $e\in\allseries$ with $x_1\not\in\supp(e)$.
Furthermore, this notion coincides with the classical definition when $y=F_c[u]$
has an analytic state space realization of the form given in \rref{eq:control-affine-realization}.
Specifically, in light of the identity
$
\dot{y}=F_{x_0^{-1}(c)}[u]+uF_{x_1^{-1}(c)}[u]
$
and \rref{eq:cF-relative-degree-decomposition}, it follows that if $c$ has relative degree $r$ then on some interval $[0,T)$ with $T>0$
and for any sufficiently small $u\in L_1[0,T]$:
\begin{subequations}
\begin{align}
y^{(k)}&= F_{(x_0^{k})^{-1}(c)}[u] \label{eq:yrminus1-in-Fliess-operators},\;\;k=0,1,\ldots,r-1 \\
y^{(r)}&=F_{(x_0^r)^{-1}(c)}[u]+uF_{(x_0^{r-1}x_1)^{-1}(c)}[u]. \label{eq:yr-in-Fliess-operators}
\end{align}
\end{subequations}
This would imply that $L_{g_1}L_{g_0}^kh(z(t))=F_{(x_0^{k}x_1)^{-1}(c)}[u]=0$ on $[0,T)$ for
$k=0,1,\ldots,r-2$. In addition, $L_{g_1}L_{g_0}^{r-1}(z_0)=F_{(x_0^{r-1}x_1)^{-1}(c)}[u](0)=(c,x_0^{r-1}x_1)=K\neq 0$.
Hence, the realization has relative degree $r$ at $z_0$ in the classical sense. Furthermore,
if the word $x_0^{r-1}x_1\notin\supp(c)$ then
$F_{(x_0^{r-1}x_1)^{-1}(c)}[u](0)=0$, but $F_{(x_0^{r-1}x_1)^{-1}(c)}[u]$ is not identically
zero on some interval $[0,T_1)$ with $T_1>0$ unless $(x_0^{r-1}x_1)^{-1}(c)=x_1^{-1}(e)=0$ (since generating series are unique).
But this can only happen if $\supp(e)\subseteq x_0X^\ast$, which contradicts the assumption that $r$ is the largest integer such that
$\supp(c_F)\subseteq x_0^{r-1}X^\ast$.
Therefore, the realization can not have a well defined
relative degree in this case.
Conversely, if the realization \rref{eq:control-affine-realization} has relative degree $r$ at $z_0$ then
for $k=0,1,\ldots,r-2$ it follows that
$L_{g_1}L_{g_0}^kh(z(t))=F_{(x_0^{k}x_1)^{-1}(c)}[u]=0$ on a neighborhood of $z_0$. So necessarily
$(x_0^{k}x_1)^{-1}(c)=0$ for $k=0,1,\ldots,r-2$, or equivalently, no word in the support of $c_F$ can have a
prefix of the form $x_0^kx_1$ with $k<r-1$. Therefore, $\supp(c_F)\subseteq x_0^{r-1}X^\ast$. In addition, since
$L_{g_1}L_{g_0}^{r-1}(z_0)=F_{(x_0^{r-1}x_1)^{-1}(c)}[u](0)=(c,x_0^{r-1}x_1)\neq 0$, then $x_0^{r-1}x_1\in\supp(c_F)$, and $x_0^{r-1}$ is
the longest such prefix that all the words in $\supp(c_F)$ can share. So $c$
must have the form of \rref{eq:cF-relative-degree-decomposition}, implying that $c$ has relative degree $r$.

\begex
Consider the state space realization \rref{eq:control-affine-realization}, where
\begdi
g_0(z)=
\left[
\begin{array}{c}
z_2^2+z_2z_3+z_3 \\
z_1^5+z_3 \\
z_1^2
\end{array}
\right],\;\;
g_1(z)=
\left[
\begin{array}{c}
0 \\
0 \\
1
\end{array}
\right],\;\;
h(z)=z_1.
\enddi
Observe that
$L_{g_1}h(z)=0$ and $L_{g_1}L_{g_0}h(z)=z_2+1$,
so that the realization has relative degree 2 at $z_0=[z_{01}\;z_{02}\;z_{03}]^T$ if $z_{02}\neq -1$ and is undefined otherwise.
So, for example, the generating series for $y=F_c[u]$ when $z=[0\;1\;1]^T$ is
\begdi
c=3 x_0 + 3 x_0^2 + 2 x_0 x_1 + 2 x_0^3 + 4 x_0^2 x_1 +
  x_0 x_1 x_0 + 36 x_0^4 + 4 x_0^3 x_1 +
 2 x_0^2 x_1 x_0 +\cdots,
\enddi
which has relative degree 2 since $\supp(c_F)\subset x_0X^\ast$, $\supp(c_F)\not\subseteq x_0^{r-1}X^\ast$ for $r>2$, and $x_0x_1\in \supp(c)$.
But when $z=[0\;-1\;1]^T$, the generating series becomes
\begdi
c=x_0 - x_0^2 + 2 x_0^3 + x_0 x_1 x_0 +
 4 x_0^3 x_1 + 2 x_0^2 x_1 x_0 +
 2 x_0^2 x_1^2 + x_0 x_1 x_0 x_1 +
 6 x_0^5 + 4 x_0^3 x_1^2 +\cdots,
\enddi
which
fails to have a well defined relative degree since $\supp(c_F)\subset x_0X^\ast$ and $\supp(c_F)\not\subseteq x_0^{r-1}X^\ast$ for $r>2$, but $x_0x_1\not\in \supp(c)$.
\endex

The following theorem describes feedback input-output linearization without any requirement that the plant has
a state space realization.

\begth
Suppose $c\in\allseriesLC$ has relative degree $r$. Then $y=F_{c}[u]$ is feedback input-output linearized by $u=F_{e_\delta^{-1}}[v]$, where
$e_\delta=((x_0^{r-1}x_1)^{-1}(c),(x_0^r)^{-1}(c))$ and provided that $e_\delta^{-1}\in\allseriesdeltaLC$,
giving the closed-loop system $y=v^{(-r)}$.
\endth

\begpr
It follows from \rref{eq:yr-in-Fliess-operators} that $y^{(r)}=F_{e_\delta}[u]$. Therefore, setting $u=F_{e_\delta^{-1}}[v]$ gives
$
y^{(r)}=F_{e_\delta}\circ F_{e_\delta^{-1}}[v]=F_{e_\delta\circ e_\delta^{-1}}[v]=F_{(1,0)}[v]=v,
$
which yields the desired result.
\endpr

\begex
Consider the polynomial $c=x_1+x_1^2$, which has relative degree 1. Setting $e_\delta=(x_1^{-1}(c),x_0^{-1}(c))=(1+x_1,0)$,
the inverse of $e_\delta$ can be computed by the same computational approach as in Example~\ref{ex:multiplicative-feedback-formula-computed}:
\begin{align*}
(1+x_1,0)^{-1}&=([\mathbf{1}+x_1Sb_{x_1}+x_1^2 Sb_{x_1^2}+x_1^3 Sb_{x_1^3}+\cdots](1+x_1),0) \\
&=(1-x_1+3x_1^2-15x_1^3+\cdots,0).
\end{align*}
Therefore, a direct calculation gives $c\modcomp (1+x_1,0)^{-1}=x_1$ as desired.
This can also be verified by observing that $c=x_1\modcomp (1+x_1,0)=\phi_d(x_1)(1)+\phi_d(x_1^2)(1)=x_1+x_1^2$.
In addition, $c\modcomp e_{\delta}^{-1}$ clearly has the same relative degree as $c$, which in the next section is shown to hold in general.
\endex

\subsection{Feedback Invariants}

It is shown next that relative degree as defined in Definition~\ref{de:relative-degree-c} is invariant under the action of
the output affine feedback transformation group. In addition, this action is free when restricted to all series having well defined relative degree.

\begth \label{th:relative-degree-invariance}
A series $c$ has relative degree $r$ if and only if it is on the
orbit of $c_N+x_0^{r-1}x_1$ under $\allseriesdeltanp$.
\endth

\begpr
If $c$ has well defined relative degree $r$ then it can be decomposed
as in \rref{eq:cF-relative-degree-decomposition}, where
without
loss of generality $e=x_0e_0+x_1e_1$ with $e_1$ proper. Then, setting
$e_\delta:=(K+e_1,e_0)\in\allseriesdeltanp$ (since $K+e_1$ is nonproper), it follows from \rref{eq:phi-definition} that
\begdi
c=c_N+x_0^{r-1}x_1(K+e_1)+x_0^re_0=c_N+\phi_e(x_0^{r-1}x_1)(1)
=(c_N+x_0^{r-1}x_1)\modcomp e_\delta.
\enddi
In which case,
$c\modcomp e_\delta^{-1}=c_N+x_0^{r-1}x_1$, or equivalently, $c$ is on the orbit of $c_N+x_0^{r-1}x_1$
under $\allseriesdeltanp$.
The converse holds since
all the steps above are reversible.
\endpr

When the series $e_\delta$ and $e_\delta^{-1}$ above
are both locally convergent then
it follows directly that if
$u=F_{e_{\delta}^{-1}}[v]$ then
$y=F_c[u]=F_c\circ F_{e_\delta^{-1}}[v]=F_{c\modcomp e_\delta^{-1}}[v]=F_{c_N+x_0^{r-1}x_1}[v]$,
as expected.

\begth
The transformation group $\allseriesdeltanp$ acts freely on the subset of
$\allseries$ having well defined relative degree.
\endth
\begpr
Assume $c$ has relative degree $r$. Without loss of generality let $c_N=0$. Then
there exists an $e_\delta\in\allseriesdeltanp$
such that $c\modcomp e_\delta^{-1}=x_0^{r-1}x_1$. So if $c\modcomp d_\delta=c$
for some $d_\delta\in\allseriesdeltanp$, then
it follows immediately that
$(c\modcomp d_\delta)\modcomp e_\delta^{-1}=c\modcomp e_\delta^{-1}$ and
$(c\modcomp e_\delta^{-1})\modcomp d_\delta^e=c\modcomp e_\delta^{-1}$,
where $d_\delta^e$ corresponds to the conjugate action
$e_\delta\circ d_\delta\circ e_\delta^{-1}$.
In which case,
$x_0^{r-1}x_1\modcomp d_\delta^e=x_0^{r-1}x_1$, or equivalently,
$x_0^{r-1}x_1d_L^e+x_0^rd_R^e=x_0^{r-1}x_1$,
and therefore, $d_\delta^e:=(d_L^e,d_R^e)=(1,0)$, the identity element of $\allseriesdeltanp$.
Thus, $e_\delta\circ d_\delta \circ e_\delta^{-1}=(1,0)$, which gives the
desired conclusion that $d_\delta=(1,0)$.
\endpr

\section{The Lie Algebra ${\bf Lie}(\mbf{\allseriesdeltanp})$}
\label{sec:FPS-Lie-group}

The main goal in this section is to describe the Lie algebra ${\rm Lie}(\allseriesdeltanp)$ of infinitesimal characters associated with the group $\allseriesdeltanp$ (see Proposition \ref{prop:exp}) in terms of a pre-Lie product.
The latter is considered the most elementary combinatorial building block since the group can be fully
reconstructed solely in terms of its pre-Lie product \cite{Foissy_EJM15,Oudom-Guin_08}.
The reader is referred to \cite{Cartier_10,Manchon_11} for details on pre-Lie algebras.
There are at least two ways to approach this problem. One way is to view $\allseriesdeltanp$ as an
affine group scheme, which in turn always corresponds to a Hopf algebra \cite[p.~9]{Waterhouse_79}, in this case the Hopf algebra
$H$ as described in
Theorem~\ref{th:hopf-algebra-for-coordinate-functions}.
In this setting, an $\re$-derivation $D$ is any linear mapping taking $H$ into an $H$-module $M$ satisfying
$D(h_1h_2)=D(h_1)h_2+h_1D(h_2)$ for all $h_1,h_2\in H$ with $D(k)=0$ when $k\in\re$. It is described in terms of the {\em universal derivation}
on $H$
$ \label{eq:universal-derivation}
{\mathbf d}=(\id\otimes\pi)\circ \Delta,
$
where $\pi:H \rightarrow H_+$ with $H_+=\bigoplus_{n > 0} H_n$ being the augmentation ideal $\mathsf{Ker}(\varepsilon)$ (recall the augmentation map
$\varepsilon:H\rightarrow \re$ of a Hopf algebra).
A linear operator $T:H\rightarrow H$
is said to be {\em left-invariant} when $\Delta\circ T=(\id\otimes T)\circ \Delta$.
In which case, the Lie algebra of the group represented by $H$ is the
$\re$-vector space of all left-invariant
derivations on $H$. This is the approach taken in \cite{Foissy_CiA15} for the case of output feedback. The
Lie algebra ${\rm Lie}(\allseriesdeltanp)$
is then determined by dualizing this approach using the fact that $a_{\eta}(\nu)={\mathds 1}_{\eta\nu}$ for all
$a_\eta\in H$ and $\nu\in X^\ast$.

A second approach is to treat $\allseriesdeltanp$ as though it
were a Lie group and to show that the Lie algebra associated
with the set of all left-invariant vector fields at the identity element can be described in terms of a pre-Lie product.
Of course, this in no way implies that $\allseriesdeltanp$ is a Lie group, but it provides a complementary point of
view to the purely algebraic first approach.
It is convenient in this context to write $c_\delta=(c_L,c_R)\in\allseriesdeltanp$ as $c_\delta=\delta c_L+c_R$,
so that the symbol $\delta$ is treated more like a letter in $X$.
The first task is to describe a left-invariant vector field on $\allseriesdeltanp$.
The left translation of $d_\delta$
by ${c_\delta}$
is $c_{\delta}\circ d_{\delta}=\delta[(c_L\modcomp d_\delta)\shuffle d_L]+[(c_L\modcomp d_\delta)\shuffle d_R+c_R\modcomp d_\delta]$,
which is denoted by
$(c_\delta\circ):\allseriesdeltanp\rightarrow\allseriesdeltanp:d_\delta\mapsto c_\delta\circ d_\delta$.
Since composition is left linear, there is no loss of generality
in setting $c_\delta=\xi_\delta:=\delta \xi_L+\xi_R$, where $\xi_L,\xi_R\in X^\ast$.
The differential of $(\xi_{\delta}\circ):\allseriesdeltanp\rightarrow\allseriesdeltanp$
at the identity element $\delta$ is the linear map $(\xi_\delta\circ)_\ast: T_\delta\allseriesdeltanp\rightarrow T_{\xi_\delta}\allseriesdeltanp$.
Consider for some $\epsilon>0$ a differentiable
path $\gamma:(-\epsilon,\epsilon)\rightarrow \allseriesdeltanp:t\mapsto d_\delta(t)$ such that $d_\delta(0)=\delta$.
Define the velocity vector at $t=0$ as the series in $\allseriesdelta$ of the form
$
v_\delta=\dot{d}_\delta(0)=\delta\dot{d}_L(0)+\dot{d}_R(0)=\delta v_L+v_R,
$
where $v_L\in\allseries$ is proper.
Then specifically the differential of $\xi_\delta\circ$ at $\delta$ in the direction of $v_\delta$ is
\begin{align*}
(\xi_\delta\circ)_\ast(v_\delta) &=\left.\frac{d}{dt} \xi_\delta \circ d_\delta(t)\right|_{t=0} \\
&= \left.\frac{d}{dt}
\delta[(\xi_L\modcomp d_\delta(t))\shuffle d_L(t)]+
(\xi_L\modcomp d_\delta(t))\shuffle d_R(t)+\xi_R\modcomp d_\delta(t)
\right|_{t=0} \\
&=
\delta\left[\left.\frac{d}{dt}\xi_L\modcomp d_\delta(t)\right|_{t=0}+\xi_L\shuffle v_L\right]+
\xi_L\shuffle v_R+\left.\frac{d}{dt}\xi_R\modcomp d_\delta(t)\right|_{t=0}.
\end{align*}
The time derivative of the product $\xi\modcomp d_{\delta}$ is computed inductively. It is clearly zero when
$\xi=\emptyset$. Otherwise, using Lemma~\ref{le:mixed-product-properties} (4),
\begin{align*}
\left.\frac{d}{dt} (x_0\xi)\modcomp d_{\delta}(t) \right|_{t=0}&=x_0\left.\frac{d}{dt} \xi\modcomp d_{\delta}(t) \right|_{t=0} \\
\left.\frac{d}{dt} (x_1\xi)\modcomp d_{\delta}(t)\right|_{t=0}
&=\left. x_1\frac{d}{dt} (d_L(t)\shuffle (\xi\modcomp d_\delta(t))) + x_0(d_R(t)\shuffle (\xi \modcomp d_\delta(t)))\right|_{t=0} \\
&=x_1\left(v_L\shuffle \xi+\left.\frac{d}{dt} \xi\modcomp d_\delta(t) \right|_{t=0}\right)+x_0(v_R\shuffle \xi).
\end{align*}
Therefore, $\left.\frac{d}{dt} \xi\modcomp d_\delta(t)\right|_{t=0}=\xi\bullet v_\delta$, where $\emptyset\bullet v_{\delta}=0$ and
\begeq \label{eq:bullet-def-for-xi}
(x_0\xi)\bullet v_\delta=x_0(\xi\bullet v_\delta),\;\;
(x_1\xi)\bullet v_\delta=x_1(v_L\shuffle \xi+\xi\bullet v_\delta)+x_0(v_R\shuffle\xi).
\endeq
So the differential in question is
\begdi
(\xi_\delta \circ)_\ast(v_\delta)=\delta[\xi_L\bullet v_\delta+\xi_L\shuffle v_L]+\xi_L\shuffle v_R+\xi_R\bullet v_{\delta}
=(\delta \xi_L)\bullet v_{\delta}+\xi_R\bullet v_{\delta}=\xi_\delta\bullet v_\delta,
\enddi
where the definition in \rref{eq:bullet-def-for-xi} is extended to treat the {\em letter} $\delta$ as
$(\delta \xi)\bullet v_\delta=\delta(v_L\shuffle \xi+\xi\bullet v_\delta)+(v_R\shuffle \xi)$.
In which case, the left-invariant vector field satisfying $\chi^{v_{\delta}}_\delta=v_\delta$ is
\begeq
\chi^{v_\delta}:\allseriesdeltanp\rightarrow {\rm T}\allseriesdeltanp:c_\delta \mapsto c_\delta\bullet v_{\delta}. \label{eq:left-invariant-vector-field}
\endeq
Now $H$ is free as a commutative algebra, that is, $H=S(V)$, where $S(V)$ is the symmetric algebra over the vector space $V$. Also
the mixed composition product is left linear, and in light of \rref{eq:grading-mixed-coproduct} the corresponding
coproduct satisfies $\tilde{\Delta}V\subseteq V\otimes H$.
Hence, $H$ is a free commutative, right-sided Hopf algebra. In which case, the graded dual of $V$, namely the vector space $V^\ast$ formed from the $\re$ span over
words in $X^\ast$ (the duals of the $a_\eta$'s) and nonempty words with prefix $\delta$, denoted by $\delta X^+$, (the duals of $b_\eta$'s, recall $b_\emptyset\sim 1$)
must have a right pre-Lie product \cite{Cartier_10,Manchon_11,Menous-Patras_16,Loday-Ronco_10}.

\begle
The vector space $V^\ast$ with product $\bullet$
is a right pre-Lie algebra, i.e.,
\begeq
(v_\delta^1\bullet v_\delta^2)\bullet v_\delta^3 - v_\delta^1\bullet (v_\delta^2\bullet v_\delta^3)=(v_\delta^1\bullet v_\delta^3)\bullet v_\delta^2 - v_\delta^1\bullet(v_\delta^3\bullet v_\delta^2). \label{eq:pre-Lie-identity}
\endeq
\endle
\begpr
The identity can be verified directly using the distributive
property $(\eta\shuffle \xi)\bullet v_\delta=(\eta\bullet v_\delta)\shuffle \xi+\eta\shuffle (\xi\bullet v_\delta)$,
which can be proved by induction on the sum of the lengths of $\eta,\xi \in X^\ast$.
\endpr

\begex
Consider \rref{eq:pre-Lie-identity} where $v^1_\delta=\delta x_1$, $v^2_\delta=x_1$ and $v^3_\delta=\delta x_0$.
Then $\delta x_1\bullet x_1=\delta x_0x_1+2x_1^2$,
$x_1\bullet \delta x_0= x_1x_0$,
$\delta x_1 \bullet \delta x_0=\delta(2x_1x_0+x_0x_1)$,
$\delta x_0\bullet x_1=x_0x_1+x_1x_0$,
and both sides of \rref{eq:pre-Lie-identity} equal
$\delta(2x_0^2x_1+x_0x_1x_0)+2x_1^2x_0+x_1x_0x_1$.
\endex

Finally, a pre-Lie product is Lie admissible, that is, it always defines a proper Lie bracket $[\cdot,\cdot]_\bullet$ under antisymmetrization. 
In the present case, this is equivalent to the Lie algebra induced by \rref{eq:left-invariant-vector-field} since
\begdi \label{eq:composition-LieBracket}
[v^1_\delta,v^2_\delta]=\left.\left[\chi^{v_\delta^1},\chi^{v_\delta^2}\right]\right|_{\delta}
=\left.\partial \chi^{v_\delta^1}(c_\delta\bullet v_\delta^2)-
\partial \chi^{v_\delta^2}(c_\delta\bullet v_\delta^1)\right|_{c_\delta=\delta}
=v_\delta^2\bullet v_\delta^1-v_\delta^1\bullet v_\delta^2.
\enddi
The next theorem follows directly from what is known for free commutative right- (or left-) sided polynomial Hopf algebras (see \cite[Theorem 5.8]{Loday-Ronco_10},
\cite[Proposition 5]{Menous-Patras_16}, and \cite[Section 3]{Oudom-Guin_08}.
\begth
The Lie algebra ${\rm Lie}(\allseriesdeltanp)$ of infinitesimal characters associated with the group $\allseriesdeltanp$ is isomorphic as a Lie algebra to $(V^\ast,[\cdot,\cdot]_\bullet)$.
\endth

\begex
In the special case where $c_\delta=\delta+c_R$ and $d_\delta=\delta+d_R$, the corresponding subspace of $T_\delta\allseriesdeltanp$ is spanned by vectors of the form $v_\delta=\delta 0+v_R$. Thus, the pre-Lie product and Lie bracket above reduce to those described in
\cite{Foissy_CiA15} and \cite{Gray-etal_CDC14}, respectively.
\endex

\section*{Acknowledgments}

The first author was supported by grant SEV-2011-0087 from the Severo Ochoa Excellence Program at the Instituto de Ciencias Matem\'{a}ticas in Madrid, Spain.
This research was also supported by a grant from the BBVA Foundation.

\end{document}